\documentclass[11pt]{article}
\usepackage{amsmath}
\usepackage{amsthm}
\let\comment\iffalse

\usepackage{amsfonts}

\newif\ifconcmath
\concmathtrue

\ifconcmath
  \usepackage[mathbf,mathcal]{euler}

  \AtBeginDocument{
    \DeclareSymbolFont{operators}   {OT1}{ccr} {m}{n}
    \SetSymbolFont{operators} {bold}{OT1}{ccr} {bx}{n}
    \DeclareMathAlphabet{\mathbold}{OT1}{ccr}{bx}{n}
  }

  \newtheoremstyle{conc-math}
    {}{}{\usefont{OT1}{ccr}{m}{n}}{}{\bfseries}{.}{.5em}{}
\else
  \AtBeginDocument{
    \DeclareMathAlphabet{\mathbf}{OML}{cmm}{b}{it}
    \DeclareMathAlphabet{\mathbold}{OT1}{cmr}{bx}{n}
  }
\fi

\newtheoremstyle{note}
  {}{}{\small}{}{\small\bfseries}{:}{ }{}

\theoremstyle{note}

\ifconcmath
  \theoremstyle{conc-math}
\else
  \theoremstyle{plain}
\fi

\newtheorem{theorem}{Theorem}[section]
\newtheorem{corollary}[theorem]{Corollary}
\newtheorem{lemma}[theorem]{Lemma}
\newtheorem{proposition}[theorem]{Proposition}

\newtheorem{conjecture}{Conjecture}

\newtheorem{definition}[theorem]{Definition}

\newenvironment{proof-of-theorem}[1]{\noindent{\bf Proof of Theorem #1}\hspace*{1em}}{\qed\bigskip}

\def\defeq{\stackrel{\mathrm{def}}{=}}

\def\Reals#1{\mathbb{R}^{#1}}

\def\vs#1#2#3{#1_{#2},\dots,#1_{#3}}

\def\abs#1{\left| #1 \right|}

\def\pnorm#1#2{\left\| #2 \right\|_{#1}}

\def\norm#1{\left\| #1 \right\|}

\def\maxnorm#1{\| #1 \|_{\max} }

\def\form#1#2{#1^{T} #2}

\def\bigoh#1{\mathcal{O}\left(#1\right)}

\def\diff#1{\, d #1 \,}

\def\setof#1{\left\{{\let\st\colon #1 }\right\}}

\def\orig#1{\bar{#1}}

\def\origin{{\mbox{\boldmath $0$}}}

\DeclareMathOperator*\erfc{erfc}
\def\prob#1#2{\Pr_{#1}\left[ #2 \right]}
\DeclareMathOperator\E{E}
\def\expec#1#2{\E_{#1}\left[ #2 \right]}

\def\laplace#1{\mathop{\mathcal{L}}[#1]}

\def\aa{\mathbf{a}}
\def\AA{\mathbf{A}}
\def\bb{\mathbf{b}}
\def\BB{\mathbf{B}}

\def\CC{\mathbf{C}}
\def\DD{\mathbf{D}}
\def\ee{\mathbf{e}}

\def\GG{\mathbf{G}}

\def\LL{\mathbf{L}}
\def\II{\mathbf{I}}

\def\QQ{\mathbf{Q}}
\def\SS{\mathbf{S}}
\def\tt{\mathbf{t}}
\def\TT{\mathbf{T}}
\def\uu{\mathbf{u}}
\def\UU{\mathbf{U}}
\def\vv{\mathbf{v}}

\def\xx{\mathbf{x}}

\def\YY{\mathbf{Y}}

\numberwithin{equation}{section}

\begin{document}

\title{Smoothed Analysis of the Condition Numbers
and  Growth Factors of Matrices}

\author{
Arvind Sankar
\thanks{Partially supported by 
NSF grant CCR-0112487}\\ 
Department of Mathematics \\ 
Massachusetts Institute of Technology
\and
Daniel A. Spielman 
\thanks{Partially supported by an Alfred P. Sloan Foundation Fellowship,
and NSF grants CCR-0112487 and CCR-0324914}\\ 
Department of Computer Science\\ 
Yale University\\
\and 
Shang-Hua Teng 
\thanks{
Partially supported by an Alfred P. Sloan Foundation Fellowship,
 NSF grant CCR-9972532, and NSF grants CCR-0311430
and ITR CCR-0325630.}\\
Department of Computer Science\\
Boston University and\\
Akamai Technologies Inc.\\
}


\maketitle

\begin{abstract}
{
Let $\orig{\AA}$ be an arbitrary matrix and let $\AA$ be a slight random
  perturbation of $\orig{\AA}$.
We prove that it is unlikely that $\AA$ has large condition number.
Using this result, we prove it is unlikely that 
  $\AA$ has large growth factor under Gaussian elimination
  without pivoting.
By combining these results, we show that the smoothed precision
  necessary to solve $\AA\xx =\bb$, for any $\bb$,
  using Gaussian elimination without
  pivoting is logarithmic.
Moreover, when $\orig{\AA}$ is an all-zero square matrix,
  our results significantly 
  improve the average-case analysis of Gaussian elimination
  without pivoting performed by Yeung and Chan (SIAM J. Matrix
  Anal. Appl., 1997).
}
\end{abstract}

\thispagestyle{empty}

\newpage
\setcounter{page}{1}

\section{Introduction}\label{sec:intro}

Spielman and Teng~\cite{SpielmanTengSimplex}, 
  introduced the smoothed analysis
  of algorithms to explain the success of algorithms
  and heuristics that could not be well understood through traditional
  worst-case and average-case analyses.
Smoothed analysis is a hybrid of worst-case and average-case analyses
  in which one measures the maximum over inputs of the expected
  value of a measure of the performance of an algorithm
  on slight random perturbations of that input.
For example, the smoothed complexity of an algorithm is the maximum over 
  its inputs of the expected running time of the
  algorithm under slight perturbations of that input.
If an algorithm has low smoothed complexity and its
  inputs are subject to noise, then it is unlikely
  that one will encounter an input on which the
  algorithm performs poorly.
(See also the Smoothed Analysis Homepage~\cite{SmoothedAnalysisHomepage})

Smoothed analysis is motivated by the existence of algorithms and
  heuristics that are known to work well in practice, but
  which are known to have poor worst-case performance.
Average-case analysis was introduced in an attempt
  to explain the success of such heuristics.
However, average-case analyses are often unsatisfying as the
  random inputs they consider may bare little resemblance to the 
  inputs actually encountered in practice.
Smoothed analysis attempts to overcome this objection by 
  proving a bound that holds in every neighborhood of inputs.

In this paper, we prove 
  that perturbations of arbitrary matrices are unlikely to
  have large condition numbers or large growth factors
  under Gaussian Elimination without pivoting.
As a consequence, we conclude that the smoothed precision
  necessary for Gaussian elimination is logarithmic.
We obtain similar results for perturbations that affect only
  the non-zero and diagonal entries of symmetric matrices.
We hope that these results will be a first step toward a smoothed
  analysis of Gaussian elimination with partial pivoting---an
  algorithm that is widely used in practice but known to have poor
  worst-case performance.
  
In the rest of this section, we recall the definitions of
  the condition numbers and growth factors of matrices, and review prior
  work on their average-case analysis.
In Section~\ref{sec:cond}, we perform a smoothed
  analysis of the condition number of a matrix.
In Section~\ref{sec:growth}, we use the results of
  Section~\ref{sec:cond} to obtain a smoothed analysis of
  the growth factors of Gaussian elimination without pivoting.
In Section~\ref{sec:smoothed}, we combine these results to 
  obtain a smoothed bound on the precision needed by Gaussian elimination
  without pivoting.
Definitions of zero-preserving perturbations and our results
  on perturbations that only affect the non-zero and diagonal
  entries of symmetric matrices appear in Section~\ref{sec:symmetric}.
In the conclusion section, 
  we explain how our results may be extended to larger families of
  perturbations,
  present some counter-examples, and
  suggest future directions for research.
Other conjectures and open questions appear in the body of the paper.

The analysis in this paper requires many results from probability.
Where reasonable, these have been deferred to the appendix.

\subsection{Condition numbers and growth factors}\label{ssec:cg}
We use the standard notation for the 1, 2 and $\infty$-norms of matrices
  and column vectors, and define
\begin{align*}
\maxnorm{\AA} & = \max_{i,j} \abs{\AA_{i,j}}.
\end{align*}

\begin{definition}[Condition Number]
For a square matrix $\AA$, the condition number of $\AA$ is defined
  by
\[
  \kappa (\AA) = \pnorm{2}{\AA} \pnorm{2}{\AA^{-1}}.
\]
\end{definition}

The condition number measures how much the solution to a system
  $\AA \xx = \bb$ changes as one makes slight changes to $\AA$ and $\bb$.
A consequence is that if
  ones solves the linear system using fewer
  than $\log (\kappa (\AA))$ bits of precision, one is likely
  to obtain a result far from a solution.
For more information on the condition number of a matrix,
  we refer the reader to one
  of~\cite{GolubVanLoan,TrefethenBau,DemmelBook}.

The simplest and most often implemented method of solving
  linear systems is Gaussian elimination.
Natural implementations of Gaussian elimination
  use $\bigoh {n^{3}}$ arithmetic operations
  to solve a system of $n$ linear equations
  in $n$ variables.
If the coefficients of these equations
  are specified using $b$ bits, 
  in the worst case it suffices to perform
  the elimination using $O(bn)$ bits of 
  precision~\cite{GLS}.
This high precision may be necessary because
  the elimination may produce large intermediate
  entries \cite{TrefethenBau}.
However, 
  in practice one usually obtains accurate
  answers using much less precision.
In fact, it is rare to find an implementation of Gaussian
  elimination that uses anything more than double precision,
  and high-precision solvers are rarely used or needed
  in practice~\cite{TrefethenBau,TrefethenSchreiber} 
  (for example, LAPACK uses 64 bits~\cite{LAPACK}).
One of the main results of this paper is 
  that $\bigoh {b + \log n}$ bits of precision usually suffice
  for Gaussian elimination in the smoothed analysis framework.

Since Wilkinson's seminal work~\cite{Wilkinson}, it has been
  understood that it suffices to carry out Gaussian elimination
  with $b + \log_{2}(5 n \kappa (\AA) \norm{\LL}_{\infty } 
  \norm{\UU}_{\infty } / \norm{\AA}_{\infty } + 3)$ 
  bits of accuracy
  to obtain a solution that is accurate to $b$ bits.
In this formula, $\LL$ and $\UU$ are the LU-decomposition of $\AA$; that is,
  $\UU$ is the upper-triangular matrix and 
  $\LL$ is the lower-triangular matrix with $1$s on the diagonal
  for which $\AA = \LL\UU$.

\subsection{Prior work}\label{ssec:prior}
The average-case behaviors of the condition numbers and growth factors of
  matrices have been studied both analytically and experimentally.
In his paper, 
``The probability that a numerical analysis problem is difficult'',  
  Demmel~\cite{DemmelProb} proved that it is unlikely that
  a Gaussian random matrix centered at the origin
  has large condition number.
Demmel's bounds on the condition number were
   improved by Edelman~\cite{Edelman}.

Average-case analysis of growth factors began with the
  experimental work of Trefethen and Schreiber~\cite{TrefethenSchreiber}, 
 who found that Gaussian random
  matrices rarely have large growth factors
  under partial or full pivoting.

\begin{definition}[Gaussian Matrix]
A matrix $\GG$ is a {\em Gaussian random matrix} of variance $\sigma^{2}$
  if each entry of $\GG$ is an independent univariate Gaussian variable
  with mean $0$ and standard deviation $\sigma $.
\end{definition}

Yeung and Chan~\cite{YeungChan} study the growth factors
  of Gaussian elimination without pivoting 
  on Gaussian random matrices of variance $1$.
They define $\rho _{\UU}$ and $\rho _{\LL}$ by
\begin{eqnarray*}
   \rho _{\UU} (\AA) & = & \norm{\UU}_{\infty } / \norm{\AA}_{\infty },
  \mbox{ and}\\
   \rho _{\LL} (\AA) & = & \norm{\LL}_{\infty },
\end{eqnarray*}
where $ \AA =\LL\UU$ is the LU-factorization of $\AA$ obtained without pivoting.
They prove
\begin{theorem}[Yeung-Chan]\label{thm:yc}
There exist constants $c > 0$ and $0 < b < 1$ such that if
  $\GG$ is an~$n \times n$ Gaussian random matrix of variance $1$ and
  $\GG = \LL\UU$ is the LU-factorization of $\GG$, then
\begin{eqnarray*}
  \prob{}{\rho _{\LL} (\GG) > x} 
& \leq &
  \frac{ c n^{3}}{x}, \mbox{ and }\\
 \prob{}{\rho _{\UU} (\GG) > x}
& \leq &
   \min \left(\frac{c n^{7/2}}{x}, \frac{1}{n} \right)
   + \frac{c n^{5/2}}{x} + b^{n}.
\end{eqnarray*}
\end{theorem}

As it is generally believed that partial pivoting is better than
  no pivoting, their result provides some intuition for the experimental
  results of Trefethen and Schreiber demonstrating that 
  random matrices rarely have large growth factors under
  partial pivoting.
However, we note that it is difficult to make this intuition rigorous
  as there are matrices $\AA$ for which no pivoting has 
  $\maxnorm{\LL} \maxnorm{\UU} / \maxnorm{\AA} = 2$ while partial pivoting
  has growth factor $2^{n-1}$. (See also~\cite{Higham})

The running times of many numerical algorithms
  depend on the condition numbers of their inputs.
For example, the number of iterations taken by the 
  method of conjugate gradients can be bounded in terms of
  the square root of the condition number.
Similarly, the running times of interior-point methods can be bounded
  in terms of condition numbers~\cite{RenegarCond}.
Blum~\cite{BlumSantaFe} suggested that
  a complexity theory of numerical algorithms should be parameterized
  by the condition number of an input in addition to the input size.
Smale~\cite{SmaleComplexity} proposed
  a complexity theory of numerical algorithms in which one:
\begin{enumerate}
\item [1.] proves a bound on the running time of an algorithm solving
  a problem in terms of its condition number, and then
\item [2.] proves that it is unlikely that a random problem
  instance has large condition number.
\end{enumerate}
This program is analogous to the average-case complexity
  of Theoretical Computer Science.

\subsection{Our results}\label{ssec:smooth}
To better model the inputs that occur in practice, 
  we propose replacing step 2 of Smale's program with
\begin{enumerate}
\item [2$'$.] prove that for every input instance it is unlikely
  that a slight random perturbation of that instance has
  large condition number.
\end{enumerate}
That is, we propose to bound the smoothed value of the condition number.
Our first result in this program is presented in Section~\ref{sec:cond},
  where we improve upon Demmel's~\cite{DemmelProb} and Edelman's~\cite{Edelman}
  average-case results to show that a slight Gaussian perturbation of an arbitrary
  matrix is unlikely to have large condition number.

\begin{definition}[Gaussian Perturbation]
Let $\orig{\AA}$ be an arbitrary $n \times n$~matrix. 
The matrix $\AA$ is a \emph{Gaussian perturbation} of $\orig{\AA}$
  of variance $\sigma^{2}$ if $\AA$ can be written as $\AA = \orig{\AA}+\GG$, where
  $\GG$ is a Gaussian random matrix of variance $\sigma^{2}$.
We also refer to $\AA$ as a Gaussian matrix of variance $\sigma^{2}$
  centered at $\orig{\AA}$.
\end{definition}

In our smoothed analysis of the condition number, we consider
  an arbitrary $n \times n$~matrix $\orig{\AA}$ of norm at most
  $\sqrt{n}$, and we bound the probability that 
  $ \kappa (\orig{\AA} + \GG)$, the condition number 
  of its Gaussian perturbation,  is large, 
  where $\GG$ is a Gaussian random matrix of variance $\sigma ^{2} \leq 1$.
We bound this probability in terms of $\sigma $ and $n$.
In contrast with the average-case analysis of 
  Demmel and Edelman, our analysis can be interpreted as demonstrating
  that if there is a little bit of imprecision or noise in the entries
  of a matrix, then it is unlikely it is ill-conditioned.
On the other hand, Edelman~\cite{EdelmanRoulette} writes of random matrices:
\begin{quotation}
What is a mistake is to psychologically link a random
  matrix with the intuitive notion of a ``typical'' matrix
  or the vague concept of ``any old matrix.''
\end{quotation}
The reader might also be interested in recent work on the smoothed
  analysis of the condition
  numbers of linear programs~\cite{BlumDunagan,DunaganSpielmanTeng,SpielmanTengTermination}.

In Section~\ref{sec:growth}, we use results from
  Section~\ref{sec:cond} to perform a smoothed analysis
  of the growth factors of Gaussian elimination without pivoting.
If one specializes our results to perturbations of an all-zero square matrix,
  then one obtains a bound on $\rho_{\UU}$ that improves the bound
  obtained by Yeung and Chan by a factor of $n$ and which
  agrees with their experimental observations.
The result obtained for $\rho_{\LL}$ also improves the bound of
   Yeung and Chan \cite{YeungChan} by a factor of $n$.
However, while Yeung and Chan compute the density functions of the
  distribution of the elements in $\LL$ and $\UU$, such precise
  estimates are not immediately available in our model.
As a result, the techniques we develop are applicable to a 
  wide variety
  of models of perturbations beyond the Gaussian.
For example, one could use our techniques to obtain results
  of a similar nature if $\GG$ were 
  a matrix of random variables chosen uniformly in $[-1,1]$.
We comment further upon this in the conclusions section of the paper.

The less effect a perturbation has, the more meaningful the
  results of smoothed analysis are.
As many matrices encountered in practice are sparse or have structure,
  it would be best to consider perturbations that respect their
  sparsity pattern or structure.
Our first result in this direction appears in Section~\ref{sec:symmetric},
  in which we consider the condition numbers and growth factors
  of perturbations of symmetric matrices that only alter their
  non-zero and diagonal elements.
We prove results similar to those proved for dense perturbations of
  arbitrary matrices.

\section{Notation and Mathematical Preliminaries}\label{sec:note}

We use bold lower-case Roman letters such as 
  $\xx $, $\aa$, $\bb_{j}$ to denote vectors in $\Reals{?}$.
Whenever a vector, say $\aa\in\Reals{n}$ is present, its components
  will be denoted by lower-case Roman letters with subscripts,
  such as $\vs{a}{1}{n}$.
Matrices are denoted by bold upper-case Roman letters 
  such as $\AA$ and scalars are
  denoted by lower-case roman letters.
Indicator random variables and random event variables
   are denoted by upper-case Roman letters.
Random variables taking real values are 
   denoted by
   upper-case Roman letters, except when they are components
   of a random vector or matrix.

The probability of an event $A$ is written $\prob{}{A}$, and the
  expectation of a variable $X$ is written $\expec{}{X}$.
The indicator random variable for an event $A$ is written $[A]$.

We write $\ln$ to denote the natural logarithm, base $e$,
  and explicitly write the base for all other logarithms.

For integers $a \leq  b$, we let $a:b$ denote the set of integers
  $\setof{x\st a\le x\le b}$.
For a matrix $\AA$ we let
  $\AA_{a:b,c:d}$
  denote the submatrix of $\AA$ indexed by 
  rows in $a:b$ and columns in $c:d$.

We will bound many probabilities by applying the following proposition.

\begin{proposition}[Minimum $\leq$ Average $\leq$ Maximum]\label{pro:avemax}
Let $\mu (X,Y)$ be a non-negative integrable function, and
 let $X$ and $Y$ be random variables distributed according to $\mu (X,Y)$.
If $A (X,Y)$ is an event and $F (X,Y)$ is a function,
  then
\begin{eqnarray*}
\min _{X} \prob{Y}{A (X,Y) }\leq  \prob{X,Y}{A (X,Y)} &\leq & \max _{X} \prob{Y}{A (X,Y)}, and \\
\min _{X} \expec{Y}{F (X,Y)}\leq   \expec{X,Y}{F (X,Y)} &\leq & \max _{X} \expec{Y}{F (X,Y)},
\end{eqnarray*}
where in the left-hand and right-hand terms, $Y$ is distributed according to 
 the induced distribution on $\mu (X,Y)$.
\end{proposition}

We recall that a matrix $\QQ$ is an orthonormal matrix if its inverse is equal to its
  transpose, that is, $\QQ^{T}\QQ = I$.
In Section \ref{sec:cond} we will use the following proposition.
\begin{proposition}[Orthonormal Transformation of Gaussian]\label{pro:orthonormal}
Let $\orig{\AA}$ be a matrix in $\Reals{n\times n}$ and $\QQ$ be an
  orthonormal matrix in $\Reals{n\times n}$.
If $\AA$ is a Gaussian perturbation of $\orig{\AA}$ of variance
  $\sigma^{2}$,
  then $\QQ\AA$ is a Gaussian perturbation of $\QQ\orig{\AA}$ of variance $\sigma^{2}$.
\end{proposition}

We will also use the following extension of Proposition 2.17 of
  \cite{SpielmanTengSimplex}.

\begin{proposition}[Gaussian Measure of Halfspaces]\label{pro:halfspace}
Let $\tt $ be any unit vector in $\Reals{n}$ and $r$ be any real.
Let $\orig{\bb}$ be a vector in $\Reals{n}$ and $\bb$ be a Gaussian
  perturbation of $\orig{\bb}$ of variance $\sigma^{2}$. 
Then
\[
\prob{\bb}{\abs{\form{\tt}{\bb}} \leq r} \leq \frac{1}{\sqrt{2\pi}\sigma
 }\int_{t=-r}^{t=r} e^{-t^{2}/2\sigma^{2}}dt.
\]
\end{proposition}


In this paper we will use the following properties of matrix norms
  and vector norms.

\begin{proposition}[Product]\label{pro:product}
For any pair of matrices $\AA$ and $\BB$ such that $\AA\BB$ is defined, and for
  every $1 \leq p \leq \infty $, 
\[
\pnorm{p}{\AA\BB}\leq \pnorm{p}{\AA}\pnorm{p}{\BB}.
\]
\end{proposition}

\begin{proposition}[Vector Norms]\label{pro:vector12}
For any column vector $\aa$ in $\Reals{n}$,
  $\norm{\aa}_{1}/\sqrt{n} \leq \norm{\aa}_{2}\leq \norm{\aa}_{1}$.
\end{proposition}

\begin{proposition}[2-norm]\label{pro:2norm}
For any matrix $\AA$,
\[
\pnorm{2}{\AA} = \pnorm{2}{\AA^{T}},
\]
as both are equal to the largest eigenvalue of $\sqrt{\AA^{T}\AA}$.
\end{proposition}

\begin{proposition}[$\norm{\AA}_{\infty}$: the maximum absolute row sum norm]\label{pro:rowsum}For every matrix $\AA$,
\begin{align*}
\norm{\AA}_{\infty} & = 
    \max_{i} \norm{\aa_i^{T}}_{1},
\end{align*}
where  $\vs{\aa}{1}{n}$ are the rows of $\AA$.
Thus, for any submatrix $\DD$ of $\AA$, 
\[
\pnorm{\infty}{\DD} \leq \pnorm{\infty}{\AA}.
\]

\end{proposition}

\begin{proposition}[$\norm{\AA}_{1}$: the maximum absolute column sum norm]\label{pro:columnsum}For every matrix $\AA$,
\begin{align*}
\norm{\AA}_{1} & = 
    \max_{i} \norm{\aa_i}_{1},
\end{align*}
where  $\vs{\aa}{1}{n}$ are the columns of $\AA$.
Thus 
\[
\pnorm{1}{\AA} = \pnorm{\infty}{\AA^{T}}.
\]

\end{proposition}

\section{Smoothed analysis of the condition number of a matrix}\label{sec:cond}

In this section, we will prove the following 
  theorem which shows that for every matrix it is unlikely that
  a slight perturbation of that matrix has large condition number.

\begin{theorem}[Smoothed Analysis of Condition number]\label{thm:kappaSmall}
Let $\orig{\AA}$ be an $n \times n$~matrix 
  satisfying $\pnorm{2}{\orig{\AA}} \leq \sqrt{n}$,
  and let $\AA$ be a Gaussian perturbation of $\orig{\AA}$ of variance 
  $\sigma ^{2} \leq  1$.
Then, $\forall x \geq 1$,
\[
  \prob{}{\kappa (\AA) \geq x }
  \leq 
\frac{14.1 n \left(1 + \sqrt{2\ln (x) / 9 n} \right)}{x \sigma }.
\]
\end{theorem}

As bounds on the norm of a random matrix are
  standard, we focus on the norm of the inverse.
Recall that 
 $1/\pnorm{2}{\AA^{-1}} = \min _{\xx } \pnorm{2}{\AA \xx }/ \pnorm{2}{\xx }$.

The first step in the proof is to bound the probability
  that~$\pnorm{2}{\AA^{-1}\vv}$ is small for a fixed unit vector~$\vv$.
This result is also used later~(in Section~\ref{ssec:rhoU}) 
  in studying the growth  factor.
Using this result and an averaging argument, we then bound the
  probability that~$\pnorm{2}{\AA^{-1}}$ is large.

\begin{lemma}[Projection of $\AA^{-1}$]\label{lem:AinversevSmall}
Let~$\orig{\AA}$ be an arbitrary square matrix in~$\Reals{n\times n}$,
  and let $\AA$ be a Gaussian perturbation of $\orig{\AA}$ of variance~$\sigma^2$.
Let~$\vv$ be an arbitrary unit vector.
Then
  \[ \prob{}{\pnorm{2}{\AA^{-1}\vv} > x} < \sqrt\frac2\pi \frac1{x\sigma} \]
\end{lemma}
\begin{proof}
Let $\QQ$ be an orthonormal matrix such that $\QQ^{T}  \ee_{1} = \vv$.
Let $\orig{\BB} =\QQ\orig{\AA}$ and $\BB =\QQ\AA$. 
By Proposition \ref{pro:orthonormal}, $\BB$ is a Gaussian perturbation
  of $\orig{\BB}$ of variance $\sigma^{2}$.
We have
\[
\pnorm{2}{\AA^{-1}\vv } = \pnorm{2}{\AA^{-1}\QQ^{T}\ee_{1} } =
\pnorm{2}{(\QQ\AA)^{-1}\ee_{1}} = \pnorm{2}{\BB^{-1}\ee_{1}}.
\]
Thus, to prove the lemma it is sufficient to show
  \[ \prob{\BB}{\pnorm{2}{\BB^{-1}\ee_{1}} > x} < \sqrt\frac2\pi \frac1{x\sigma}. \]
We observe that
  \[ \pnorm{2}{\BB^{-1}\ee_{1}} = \pnorm{2}{(\BB^{-1})_{:,1}}, \]
  the length of the first column of~$\BB^{-1}$.
The first column of~$\BB^{-1}$, by the definition of the matrix inverse, is the
  vector that is orthogonal to every row of $\BB$ but the first 
  and that has inner product~$1$~with the first row of $\BB$.
Hence its length is the reciprocal of the length of the projection of
  the first row of $\BB$ onto the subspace orthogonal to the rest of the rows.

Let $\vs{\bb}{1}{n}$ be the rows of $B$ and $\vs{\orig{\bb}}{1}{n}$
  be the rows of $\orig{\BB}$.
Note that $\bb_{i}$ is a Gaussian perturbation of $\orig{\bb}_{i}$ of
  variance $\sigma^{2}$.
Let $\tt$ be the unit vector that is 
  orthogonal to the span of $\vs{\bb}{2}{n}$.
Then
\[
\pnorm{2}{(\BB^{-1})_{:,1}} = \abs{\frac{1}{\form{\tt}{\bb_{1}}}}.
\]
Thus,
\begin{align*}
 \prob{\BB}{\pnorm{2}{\BB^{-1}\vv} > x} 
& =
  \prob{\vs{\bb}{1}{n}}{\abs{\frac{1}{\form{\tt}{\bb_{1}}}} > x} \\
& \leq   \max_{\vs{\bb}{2}{n}} \prob{\bb_{1}}{\abs{\form{\tt}{\bb_1}}<
1/x} \\
& <
  \sqrt{\frac{2}{\pi}}
  \frac{1}{x \sigma },
\end{align*}
where the first inequality follows from Proposition \ref{pro:avemax}
  and the second inequality follows from Lemma \ref{lem:distToPlane}.
\end{proof}

\begin{theorem}[Smallest singular value]\label{thm:sigmaSmall}
Let~$\orig{\AA}$ be an arbitrary square matrix in~$\Reals{n\times n}$,
  and let $\AA$ be a Gaussian perturbation of $\orig{\AA}$ of variance~$\sigma^2$.
Then
  \[ \prob{\AA}{\pnorm{2}{\AA^{-1}} \geq x} \leq 2.35 \frac{\sqrt n}{x\sigma} \]
\end{theorem}

\begin{proof} 
Let ~$\vv$ be a uniformly distributed
  random unit vector in $\Reals{n}$.
It follows from Lemma~\ref{lem:AinversevSmall} that
  \begin{equation}\label{eq:AinversevSmall}
    \prob{\AA,\vv}{\pnorm{2}{\AA^{-1}\vv}\geq x} \leq \sqrt\frac2\pi\frac1{x\sigma}
  \end{equation}

Since $\AA$ is a Gaussian perturbation of $\orig{\AA}$, 
  with probability~$1$ there
  is a unique pair $(\uu,-\uu)$ of unit 
  vectors such that~$\pnorm{2}{\AA^{-1}\uu}=\pnorm{2}{\AA^{-1}}$.
From the inequality
\[
  \pnorm{2}{\AA^{-1}\vv} \ge \pnorm{2}{\AA^{-1}}\abs{\form{\uu }{\vv }},
\]
we know that for every $c > 0$,
\begin{align*}
  \prob{\AA,\vv }{\pnorm{2}{\AA^{-1} \vv} \geq x \sqrt{c/n}}
& \geq 
  \prob{\AA,\vv}
     {\pnorm{2}{\AA^{-1}} \geq x 
          \mbox{ and }
      \abs{\form{\uu}{\vv}} \geq \sqrt{c/n}}\\
& =
  \prob{\AA,\vv }
     {\pnorm{2}{\AA^{-1}} \geq x}
  \prob{\AA,\vv}
     {\abs{ \form{\uu}{\vv}} \geq \sqrt{c/n} \quad \Bigg| \pnorm{2}{\AA^{-1}} \geq x }\\
& =   \prob{\AA}
     {\pnorm{2}{\AA^{-1}} \geq x}
  \prob{\AA,\vv}
     {\abs{ \form{\uu}{\vv}} \geq \sqrt{c/n} \quad \Bigg| \pnorm{2}{\AA^{-1}}
\geq x }\\
& \geq \prob{\AA}
     {\pnorm{2}{\AA^{-1}} \geq x}
  \min_{\AA: \pnorm{2}{\AA^{-1}}\geq x }  \prob{\vv}
     {\abs{ \form{\uu}{\vv}} \geq \sqrt{c/n}} 
  & \text{(by Proposition~\ref{pro:avemax})} \\
&\geq \prob{\AA}
     {\pnorm{2}{\AA^{-1}} \geq x}
    \prob{G}{\abs{G} \geq \sqrt{c}},  & \text{(by Lemma~\ref{lem:randSphere})}
\end{align*}
where $G $ 
   is a Gaussian random variable
   with mean $\origin$ and variance $1$.
To prove this last inequality,
  we first note that
   that $\vv $ is a
   random unit vector and is independent from $\uu$.
Thus, in 
  a basis of $\Reals{n}$ in which $\uu $ is the first vector, $\vv $ is 
   a uniformly distributed random unit vector with the first coordinate
   equal to $\form{\uu}{\vv }$,
and so we may apply
   Lemma~\ref{lem:randSphere} to bound 
  $\prob{\vv}
     {\abs{ \form{\uu}{\vv}} \geq \sqrt{c/n}}$ from below by 
   $\prob{G}{\abs{G} \geq \sqrt{c}}$.
So,
\begin{align*}
  \prob{\AA}{\pnorm{2}{\AA^{-1}} \geq x}
 &\leq 
\frac{ \prob{\AA,\vv}{\pnorm{2}{\AA^{-1} \vv} \geq x \sqrt{c/n}} }
  {\prob{G}{\abs{G} \geq \sqrt{c}}} \\
 &\leq
  \sqrt{\frac{2}{\pi}}\frac{\sqrt{n}}{x\sigma \sqrt{c}
    \prob{G}{\abs{G} \geq \sqrt{c}}}
&  \text{(by \eqref{eq:AinversevSmall})}.\\
\end{align*}
Because this inequality is true for every $c$, we will choose a value for $c$
  that almost maximizes $\sqrt{c}\prob{G}{\abs{G}\geq \sqrt{c}}$ and which in turn
  almost minimizes the right hand side.

Choosing $c = 0.57$, and evaluating the error function numerically, we determine 
\begin{align*}
  \prob{\AA}{\pnorm{2}{\AA^{-1}} \geq x}
 &\leq
   2.35 \frac{\sqrt{n}}{x\sigma}.
\end{align*}
\end{proof}

Note that Theorem \ref{thm:sigmaSmall} gives a smoothed analogue of the
  following bound of Edelman \cite{Edelman} on Gaussian random matrices.

\begin{theorem}[Edelman]\label{thm:edelman}
Let $\GG \in \Reals{n\times n}$ be a Gaussian random matrix with variance $\sigma^{2}$, then
\[
\prob{\GG}{\pnorm{2}{\GG^{-1}}\geq x} \leq \frac{\sqrt{n}}{x\sigma }.
\]
\end{theorem}

As Gaussian random matrices can be viewed as  Gaussian random perturbations of
  the $n \times n$ all-zero square matrix, Theorem \ref{thm:sigmaSmall}
  extends Edelman's theorem to Gaussian random perturbations of an
  arbitrary matrix.
The constant $2.35$ in Theorem \ref{thm:sigmaSmall} is bigger than
  Edelman's 1 for Gaussian random matrices.
We conjecture that it is possible to reduce $2.35$ in Theorem
  \ref{thm:sigmaSmall} to  1 as well.

\begin{conjecture}[Smallest Singular Value]\label{conj:sigmaSmall}
Let~$\orig{\AA}$ be an arbitrary square matrix in~$\Reals{n\times n}$,
  and let $\AA$ be a Gaussian perturbation of $\orig{\AA}$ of variance~$\sigma^2$.
Then
  \[ \prob{\AA}{\pnorm{2}{\AA^{-1}} \geq x} \leq \frac{\sqrt n}{x\sigma} \]
\end{conjecture}

We now apply Theorem \ref{thm:sigmaSmall} to prove Theorem \ref{thm:kappaSmall}.

\begin{proof} [Proof of Theorem \ref{thm:kappaSmall}]
As observed by Davidson and Szarek~\cite[Theorem II.7]{DavidsonSzarek},
  one can apply inequality (1.4) of~\cite{LedouxTalagrand} to show
  that for all $k \geq 0$,
\[
  \prob{\AA}{\pnorm{2}{\orig{\AA} - \AA} \geq \sigma \left(2\sqrt{n} + k \right)}
    \leq  e^{-k^{2}/2}.
\]
Replacing $\sigma$ by its upper bound of $1$ and setting
  $\epsilon = e^{-k^{2}/2}$, we obtain
\[
  \prob{\AA}{\pnorm{2}{\orig{\AA} - \AA} \geq 2\sqrt{n}
               + \sqrt{2 \ln (1/\epsilon )} } \leq  \epsilon,
\]
for all $\epsilon \leq 1$.
By assumption, $\pnorm{2}{\orig{\AA}} \leq \sqrt{n}$;
  so,
\[
  \prob{\AA}{\pnorm{2}{\AA} \geq 
        3 \sqrt{n}  + \sqrt{2 \ln (1/\epsilon )}}
   \leq  \epsilon.
\]
From the result of Theorem~\ref{thm:sigmaSmall}, we have
\[
\prob{\AA}{\pnorm{2}{\AA^{-1}}\geq \frac{2.35\sqrt{n}}{\epsilon \sigma }}\leq \epsilon.
\]
Combining these two bounds,  we find
\[
 \prob{\AA}{\pnorm{2}{\AA}  \pnorm{2}{\AA^{-1}} \geq 
   \frac{7.05 n + 2.35 \sqrt{2 n \ln (1/\epsilon )}}
        {\epsilon \sigma }
   }
  \leq 2 \epsilon.
\]
So that we can express this probability in the form of 
  $ \prob{\AA}{\pnorm{2}{\AA}  \pnorm{2}{\AA^{-1}} \geq x}$, for $x\geq 1$,
we let
\begin{eqnarray}\label{eqn:substi}
  x =  \frac{7.05 n + 2.35 \sqrt{2 n \ln (1/\epsilon )}}
        {\epsilon \sigma }.
\end{eqnarray}

It follows Equation (\ref{eqn:substi}) and the assumption $\sigma \leq 1$
  that $x\epsilon \geq 1$, implying $\ln (1/\epsilon ) \leq \ln x$.
From Equation (\ref{eqn:substi}), we derive
\[
 2 \epsilon 
=
\frac{2  \left(7.05 n + 2.35 \sqrt{2 n \ln (1/\epsilon )}\right)}{x \sigma }
\leq 
\frac{2  \left(7.05 n + 2.35 \sqrt{2 n \ln x}\right)}{x \sigma }
\leq 
\frac{14.1n \left(1 + \sqrt{2\ln (x) / 9 n} \right)}{x \sigma }.
\]

Therefore, we conclude 
\[
  \prob{}{\pnorm{2}{\AA}  \pnorm{2}{\AA^{-1}} \geq x }
  \leq 
\frac{14.1 n \left(1 + \sqrt{2\ln (x) /9 n} \right)}{x \sigma }.
\]
\end{proof}

We conjecture that the $1 + \sqrt{2\ln(x) / 9 n}$ term should
  be unnecessary because those matrices for which
  $\pnorm{2}{\AA}$ is large are less likely to have $\pnorm{2}{\AA^{-1}}$ large
  as well.

\section{Growth Factor of Gaussian Elimination without
Pivoting}\label{sec:growth}

We now turn to proving a bound on the growth factor.
We will consider a matrix $\AA \in \Reals{n\times n}$
  obtained from a Gaussian perturbation of variance $\sigma^{2}$ of 
  an arbitrary matrix $\orig{\AA}$ satisfying $\pnorm{2}{\orig{\AA}}\leq 1$.
With probability~$1$, none of the diagonal 
  entries that occur during elimination will be~$0$.
So, in the spirit of Yeung and Chan \cite{YeungChan},
    we analyze the growth factor of 
    Gaussian elimination without pivoting.
When we specialize our smoothed analyses to the case
  $\orig{\AA} = 0$, we improve the bounds of Yeung and Chan (see Theorem
  \ref{thm:yc})
  by a factor of $n$.
Our improved bound on $\rho _{\UU}$ agrees with their experimental analyses.

\subsection{Growth in $\UU$}\label{ssec:rhoU}

We recall that
  \[ \rho_\UU (\AA) = \frac{\norm{\UU}_\infty}{\norm{\AA}_\infty}. \]
In this section, we give two bounds on $\rho_\UU(\AA)$.
The first will have a better dependence on $\sigma $, and second
  will have a better dependence on $n$.
It is the later bound, Theorem~\ref{thm:anotherU}, that
  agrees with the experiments of Yeung and Chan~\cite{YeungChan}
  when specialized to the average-case by setting $\orig{\AA} = 0$ and $\sigma =1$.

\subsubsection{First bound}\label{sub:firstbound}

\begin{theorem}[First bound on $\rho_\UU(\AA)$]\label{thm:U}
Let $\orig{\AA}$ be an $n \times n$~matrix satisfying
  $\pnorm{2}{\orig{\AA}} \leq 1$, and let $\AA$ be a Gaussian perturbation
  of $\orig{\AA}$ of variance $\sigma ^{2} \leq 1$.
Then,
  \[ \prob{}{\rho_\UU(\AA) > 1+x} <
   \frac{1}{\sqrt{2 \pi}}\frac{n (n+1)}{x\sigma}.
   \]
\end{theorem}
\begin{proof}
By Proposition \ref{pro:rowsum}.
\[ \rho_\UU (\AA) = \frac{\norm{\UU}_\infty}{\norm{\AA}_\infty} = 
     \max_i \frac{\pnorm{1}{(\UU_{i,:})^{T}}}{\norm{\AA}_\infty}. \]

So, we need to bound the probability that  the~$1$-norm of the vector
  defined by each row of~$\UU$ is large and 
  then apply a union bound to bound the overall probability.

Fix for now a $k$ between $2$ and $n$.
We denote the upper triangular segment of 
  the~$k$th row of~$\UU$ by~$\uu^{T} = \UU_{k,k:n}$, 
and observe that $\uu$ can be obtained from the formula:
  \begin{equation}\label{eqn:abCD}
    \uu^{T} = \aa^{T} - \bb^T \CC^{-1} \DD 
  \end{equation}
  where
  \[        \aa^{T}   = \AA_{k,k:n}
     \qquad \bb^{T} = \AA_{k,1:k-1}
     \qquad \CC     = \AA_{1:k-1,1:k-1}
     \qquad \DD     = \AA_{1:k-1,k:n}.
  \]
This expression for $\uu$ follows immediately from
\[
  \AA_{1:k,:} =
    \begin{pmatrix}
      \CC    & \DD \\
      \bb^{T} & \aa^{T}
    \end{pmatrix}
  = \begin{pmatrix}
      \LL_{1:k-1,1:k-1} & 0\\
      \LL_{k,1:k-1}     & 1
    \end{pmatrix}
    \begin{pmatrix}
      \UU_{1:k-1,1:k-1} & \UU_{1:k-1,k:n}\\
      0               & \uu^{T}
    \end{pmatrix}.
\]
From~\eqref{eqn:abCD}, we derive
  \begin{align}
    \pnorm{1}{\uu} = \norm{\aa - \left(\bb^T \CC^{-1} \DD \right)^{T}}_{1}
       &\le \norm{\aa}_1 +
          \pnorm{1}{\left(\bb^T \CC^{-1}\DD \right)^{T}}\notag \\
       &\le \pnorm{\infty}{\aa^{T}} +
         \pnorm{1}{ \left(\CC^{T} \right)^{-1}\bb} \  \norm{\DD}_{\infty }\notag
  & \text{by Propositions \ref{pro:product} and \ref{pro:columnsum}}
\\
       &\le \norm{\AA}_\infty \left(1+\pnorm{1}{ \left(\CC^{T}
\right)^{-1}\bb}\right)
  & \text{by Proposition \ref{pro:rowsum}}
\label{eqn:rhoU}
  \end{align}

We now bound the probability $\pnorm{1}{\left(\CC^{T} \right)^{-1}\bb}$ is large.
By Proposition \ref{pro:vector12}, 
  
\[
\pnorm{1}{\left(\CC^{T} \right)^{-1}\bb} \leq \sqrt{k-1}\pnorm{2}{\left(\CC^{T} \right)^{-1}\bb}.
\]

Note that $\bb$ and $\CC$ are independent of each other.  Therefore,
\begin{multline}\label{eqn:thmU1}
\prob{\bb,\CC }{\pnorm{1}{\left(\CC^{T} \right)^{-1}\bb} > x} \leq 
\prob{\bb,\CC }{\pnorm{2}{\left(\CC^{T} \right)^{-1}\bb} > x/\sqrt{k-1}}
\\
\leq  
 \sqrt{\frac{2}{\pi}} \frac{\sqrt{k-1}\sqrt{(k-1)\sigma^{2}
+1}}{x\sigma} 
  <   \sqrt{\frac{2}{\pi}}\frac{k}{x\sigma},
\end{multline}
where the second inequality follows from Lemma \ref{lem:bounding1norm}
 below and the last inequality follows from the assumption
 $\sigma^{2}\leq 1$.

We now apply a union bound over the choices of $k$ to obtain
\begin{eqnarray*}
  \prob{}{\rho_\UU(\AA) > 1+x} <
    \sum_{k=2}^{n}\sqrt{\frac{2}{\pi}}\frac{k}{x\sigma}
  \leq 
   \frac{1}{\sqrt{2\pi}}\frac{n (n+1)}{x\sigma}.
\end{eqnarray*}

\end{proof}

\begin{lemma}\label{lem:bounding1norm}
Let $\orig{\CC }$ be an arbitrary square matrix in $\Reals{d\times d}$,
  and $\CC $ be a Gaussian perturbation of $\orig{\CC }$ of variance $\sigma^{2}$.
Let $\orig{\bb}$ be a column vector in $\Reals{d}$ such that
  $\norm{\orig{\bb}}_{2}\leq 1$,
and let $\bb$ be a Gaussian perturbation of $\orig{\bb}$ of variance
  $\sigma^{2}$.
If $\bb$ and $\CC $ are independent of each other, 
 then 
\[
\prob{\bb,\CC }{\pnorm{2}{\CC ^{-1} \bb}\geq x}\leq 
 \sqrt{\frac2\pi}\frac{\sqrt{\sigma^{2}d + 1}}{x\sigma} 
\]
\end{lemma}
\begin{proof}

Let~$\hat{\bb}$ be the unit vector in the direction of~$\bb$.
By applying Lemma~\ref{lem:AinversevSmall}, we obtain for all $\bb$,
  \[ \prob{\CC }{\pnorm{2}{\CC^{-1}\bb} > x}
      = \prob{\CC }{\pnorm{2}{\CC^{-1}\hat{\bb} } > \frac{x}{\norm{\bb}_2}}
      \leq
        \sqrt{\frac2\pi}\frac1{x\sigma}\norm{\bb}_{2}.
  \]

Let $\mu (\bb)$ denote the density according to which $\bb$
  is distributed.
Then, we have
\begin{align*}
\prob{\bb,\CC }{\pnorm{2}{\CC ^{-1}\bb}>x}
 & = \int_{\bb\in\Reals{d}} \prob{\CC }{\pnorm{2}{\CC ^{-1}\bb}>x} \mu (\bb)d\bb
\\
 & \leq  \int_{\bb\in\Reals{d}} \left(\sqrt{\frac2\pi}\frac1{x\sigma}\norm{\bb}_{2} \right) \mu (\bb)d\bb\\
&  = \sqrt{\frac2\pi}\frac1{x\sigma}\expec{\bb}{\norm{\bb}_{2}}.
\end{align*}

It is known~\cite[p. 277]{StatsEncyc6} that 
  $\expec{\bb}{\pnorm{2}{\bb}^{2}}\leq \sigma^{2}d +
\pnorm{2}{\orig{\bb}}^{2}$.
As  $\expec{}{X} \leq \sqrt{\expec{}{X^{2}}}$
  for every positive random variable $X$,
  we have $\expec{\bb}{\pnorm{2}{\bb}}\leq \sqrt{\sigma^{2}d +
\pnorm{2}{\orig{\bb}}^{2}}\leq \sqrt{\sigma^{2}d + 1} $.

\end{proof}

\subsubsection{Second Bound for $\rho_\UU(\AA)$}

In this section, we establish an upper bound on 
  $\rho_\UU(\AA)$ which dominates the bound in Theorem~\ref{thm:U}
  for~$\sigma\geq n^{-3/2}$.

If we specialize the parameters in this bound to $\orig{\AA}=0$
  and $\sigma ^{2} = 1$, we 
  improve the average-case bound proved by 
  Yeung and Chan~\cite{YeungChan} (see Theorem
  \ref{thm:yc}) by a
  factor of~$n$.
Moreover, the resulting bound agrees with their experimental results.

\begin{theorem}[Second bound on $\rho_\UU(\AA)$]\label{thm:anotherU}
Let $\orig{\AA}$ be an $n \times n$~matrix satisfying
  $\pnorm{2}{\orig{\AA}} \leq 1$, and let $\AA$ be a Gaussian
  perturbation of $\orig{\AA}$ of variance $\sigma ^{2} \leq 1$.
For $n\geq 2$,
  \[
   \prob{}{\rho_\UU(\AA) > 1+x} \le 
   \sqrt{\frac2\pi}\frac{1}{x}\left(
     \frac{2}{3}n^{3/2} + \frac{n}{\sigma} +
       \frac{4}{3}\frac{\sqrt{n}}{\sigma^2}\right)
   \]
\end{theorem}
\begin{proof}
As in the proof of Theorem~\ref{thm:U}, we will separately consider the
  $k$th row of $U$ for each $2 \leq k \leq n$.
For any such $k$,
 define $\uu $, $\aa$, $\bb$, $\CC $ and $\DD$ as in the proof of Theorem
 \ref{thm:U}.

In the case when $k = n$, we may apply
   \eqref{eqn:thmU1} in the proof of Theorem \ref{thm:U}, to show
\begin{eqnarray}\label{eqn:lastrow}
 \prob{}{\frac{\pnorm{1}{\uu }}{\norm{\AA}_{\infty}} > 1+x} \le 
   \sqrt{\frac{2}{\pi}}\frac{n}{x\sigma}.
\end{eqnarray}

We now turn to the case $k \leq n-1$.
By \eqref{eqn:abCD} and Proposition \ref{pro:vector12}, we have
\[
\pnorm{1}{\uu} \leq \pnorm{1}{\aa} +
          \pnorm{1}{\left(\bb^T \CC^{-1}\DD \right)^{T}} 
  \leq \pnorm{1}{\aa} +  \sqrt{k-1}\pnorm{2}{\left(\bb^T \CC^{-1}\DD
\right)^{T}} 
 = \pnorm{1}{\aa} +  \sqrt{k-1}\pnorm{2}{\bb^T \CC^{-1}\DD}.
\]
The last equation follows from Proposition \ref{pro:2norm}.
Therefore, for all $k\leq n-1$, 
\begin{align*}
\frac{\norm{\uu}_{1}}{\norm{\AA}_{\infty}} & \leq   
\frac{\norm{\aa}_{1} + \sqrt{k-1}\pnorm{2}{\bb^T \CC^{-1}\DD}}
 {\norm{\AA}_{\infty}} \leq   1 + 
\frac{ \sqrt{k-1}\pnorm{2}{\bb^T \CC^{-1}\DD}}
 {\norm{\AA}_{\infty}}& \text{(by Proposition \ref{pro:rowsum})}.\\ 
& \leq   1 + 
  \frac{ \sqrt{k-1}\pnorm{2}{\bb^T \CC^{-1}\DD}}
 {\pnorm{1}{\left(\AA_{n,:} \right)^{T}}} & \text{(also by Proposition \ref{pro:rowsum})}.
\end{align*}
We now observe that for fixed~$\bb$ and~$\CC $,
  $(\bb^{T}\CC^{-1}) \DD$ is a Gaussian
  random row vector of variance 
  $\norm{\bb^{T} \CC^{-1}}^{2}_{2}\sigma^{2}$ 
  centered at $(\bb^{T}\CC^{-1})\orig{\DD}$, where
  $\orig{\DD}$ is the center of $\DD$.
We have $\norm{\orig{\DD}}_{2}\leq \norm{\orig{\AA}}_{2}\leq 1$,
  by the assumptions of the theorem; so,
\[
\pnorm{2}{\bb^{T}\CC^{-1}\orig{\DD}} \leq 
  \norm{\bb^T\CC^{-1}}_{2}\norm{\orig{\DD}}_{2}\leq
  \norm{\bb^T\CC^{-1}}_{2}.
\]
Thus, if we let
 $\tt^{T} = (\bb^{T}\CC^{-1}\DD)/\pnorm{2}{\bb^T\CC^{-1}}$,
 then for every fixed $\bb$ and $\CC $, $\tt$ is a Gaussian random column
 vector in $\Reals{n-k+1}$ of variance $\sigma ^{2}$ centered at a vector of
 2-norm at most $1$.
We also have
\begin{eqnarray}\label{eqn:anotherU1}
\prob{\bb,\CC ,\DD}{\pnorm{2}{\bb^{T}\CC^{-1}\DD}\geq x}
& = & 
\prob{\bb,\CC ,\tt}{\norm{\bb^T\CC^{-1}}_{2}\norm{\tt}_{2} \geq x}.
\end{eqnarray}
It follows from Lemma \ref{lem:bounding1norm} that
\[
\prob{\bb,\CC }{\norm{\bb^{T}\CC^{-1}}_{2}\geq x}\leq 
 \sqrt{\frac2\pi}\frac{\sqrt{\sigma^{2} (k-1) + 1}}{x\sigma} .
\]
Hence, we may apply Corollary \ref{cor:chiComb} to show
\begin{align}
\prob{\bb,\CC ,\tt}{\norm{\bb^T\CC^{-1}}_{2}\norm{\tt}_{2} \geq x}
 & \leq \sqrt{\frac2\pi}\frac{\sqrt{\sigma^{2} (k-1) + 1}
  \sqrt{\sigma^{2} (n-k+1)+1}}{x\sigma} \notag \\
 & \leq \sqrt{\frac{2}{\pi }}\frac{\left(1+\frac{n\sigma^2}{2}\right)} {x\sigma}\label{eqn:anotherU2}.
\end{align}

Note that  $\AA_{n,:}$ is a Gaussian perturbation of variance
   $\sigma^{2}$ of a 
   row vector in $\Reals{n}$.
As $\AA_{n,:}$ is independent of $\bb$, $\CC $ and $\DD$, we
 can apply \eqref{eqn:anotherU1}, \eqref{eqn:anotherU2} and Lemma \ref{lem:linear-expect} to show
\begin{align*}
\prob{}{\frac{ \sqrt{k-1}\norm{\bb^{T}\CC^{-1}\DD}_{2}}
 {\pnorm{1}{\left( \AA_{n,:}\right)^{T}}}\geq  x }
& \leq
   \sqrt{\frac2\pi}
  \frac{\sqrt{k-1}\left(1+\frac{n\sigma^2}{2}\right)} {x\sigma}
\expec{}{\frac{1}{\pnorm{1}{\left( \AA_{n,:}\right)^{T}}}}\\
& \leq
   \sqrt{\frac2\pi}
     \frac{\sqrt{k-1}\left(1+\frac{n\sigma^2}{2}\right)} {x\sigma}
     \frac{2}{n\sigma},
\end{align*}
by Lemma~\ref{lem:expected1normGaussian}.

Applying a union bound over the choices for $k$,
  we obtain

\begin{eqnarray*}
  \prob{}{\rho_\UU(\AA) > 1+x} & \le &
   \left( \sum_{k=2}^{n-1} 
   \sqrt{\frac2\pi}
     \frac{\sqrt{k-1}\left(1+\frac{n\sigma^2}{2}\right)} {x\sigma}
     \frac{2}{n\sigma} \right)
    + \sqrt{\frac{2}{\pi}}
            \frac{n}{x\sigma}\\
 & \leq &
   \sqrt{\frac2\pi}\frac{1}{x}
   \left(
      \frac{2}{3}
     \sqrt{n}\left(\frac{2}{\sigma^2}+n\right)
    +\frac{n}{\sigma}
   \right)
\\
 & = &
   \sqrt{\frac2\pi}\frac{1}{x}\left(
     \frac{2}{3}n^{3/2} + \frac{n}{\sigma} +
       \frac{4}{3}\frac{\sqrt{n}}{\sigma^2}
     \right),
\end{eqnarray*}
where the second inequality follows from 
\[
\sum_{k=1}^{n-2}\sqrt{k} \leq \frac{2}{3}n^{3/2}.
\]
\end{proof}

\subsection{Growth in $\LL$}\label{ssec:rhoL}
Let $\LL$ be the lower-triangular part of the LU-factorization
  of $\AA$.
We have 
\[
   \LL_{(k+1):n,k} = \AA^{(k-1)}_{(k+1):n,k} \, \Big/ \AA^{(k-1)}_{k,k} ,
\]
where we let $\AA^{(k)}$ denote the matrix remaining after the
  first $k$ columns have been eliminated.
So, $\AA^{(0)} = \AA$.

Recall $\rho_{\LL} (\AA) =\pnorm{\infty}{\LL}$, which is equal to the
  maximum absolute row sum of $\LL$ (Proposition \ref{pro:rowsum}).
We will show that it is unlikely that
  $\norm{\LL_{(k+1):n,k}}_{\infty }$ is large
  by proving that it is unlikely that
  $\norm{\AA^{(k-1)}_{(k+1):n,k}}_{\infty }$ is large while 
  $\abs{\AA^{(k-1)}_{k,k}}$ is small.

\begin{theorem}[$\rho _{\LL} (\AA)$]\label{thm:L}
Let $\orig{\AA}$ be an $n$-by-$n$ matrix for which $\pnorm{2}{\orig{\AA}} \leq 1$,
  and let $\AA$ be a Gaussian perturbation of $\orig{\AA}$ of variance
$\sigma^{2} \leq 1$.
If $n\geq 2$, then, 
\[
  \prob{}{\rho _{\LL} (\AA)
          > x}
  \leq 
 \sqrt{\frac{2}{\pi }}
  \frac{n^{2}}{x}
\left(
\frac{ \sqrt{2} }{ \sigma }
+
  \sqrt{2 \ln n} + \frac{1}{\sqrt{2\pi} \ln n}
 \right)
\]
\end{theorem}
\begin{proof}
For each $k$ between $1$ and $n-1$, we have
\begin{eqnarray*}
   \LL_{(k+1):n,k} & =  &  \frac{\AA^{(k-1)}_{(k+1):n,k}}{\AA^{(k-1)}_{k,k}} \\
  & = &
   \frac{\AA_{(k+1):n,k} -  \AA_{(k+1):n, 1:(k-1)}
     \AA_{1:(k-1), 1:(k-1)}^{-1}
     \AA_{1:(k-1), k}
   }{
     \AA_{k,k} -  \AA_{k, 1:(k-1)}
     \AA_{1:(k-1), 1:(k-1)}^{-1}
     \AA_{1:(k-1), k}
   }\\
  & = &
   \frac{
    \AA_{(k+1):n,k} -  \AA_{(k+1):n, 1:(k-1)} \vv 
   }{
    \AA_{k,k} -  \AA_{k, 1:(k-1)} \vv
   },
\end{eqnarray*}
where we let $\vv =  \AA_{1:(k-1), 1:(k-1)}^{-1} \AA_{1:(k-1), k}$.
Since $\pnorm{2}{\orig{\AA}} \leq 1$, 
  and all the terms $\AA_{(k+1):n,k}$, $\AA_{(k+1):n, 1:(k-1)}$,
  $ \AA_{k,k}$, $\AA_{k, 1:(k-1)}$ and $\vv$ are independent,
   we can apply Lemma~\ref{lem:vecRatio} to show that 

\begin{align*}
  \prob{}{\norm{\LL_{(k+1):n,k}}_{\infty } > x} 
&  \leq 
 \sqrt{\frac{2}{\pi }}
  \frac{1}{x}
\left(
\frac{ \sqrt{2} }{ \sigma }
+
  \sqrt{2 \ln (\max(n-k,2))} + \frac{1}{\sqrt{2\pi} \ln (\max(n-k,2))}
 \right) \\
& \leq 
 \sqrt{\frac{2}{\pi }}
  \frac{1}{x}
\left(
\frac{ \sqrt{2} }{ \sigma }
+
  \sqrt{2 \ln n} + \frac{1}{\sqrt{2\pi} \ln n},
 \right) 
\end{align*}
where the last inequality follows the facts  that $\sqrt{2z}+\frac{1}{\sqrt{2\pi }z}$
  is an increasing function when $z \geq \pi^{-1/3} $, and
  $\ln 2 \geq \pi^{-1/3}$.

The theorem now follows by applying a union bound over
  the $n$ choices for $k$
  and observing that 
   $\norm{\LL}_{\infty }$ is at most $n$ times the largest
  entry in $\LL$.
\end{proof}

\begin{lemma}[Vector Ratio]\label{lem:vecRatio}
Let $d$ and $n$ be positive integers.
Let $a$, $\bb$, $\xx $, and $\YY$ be Gaussian perturbations of 
 $\orig{a}\in \Reals{1}$, $\orig{\bb}\in \Reals{d}$, $\orig{\xx } \in
  \Reals{n}$, and $\orig{\YY}\in \Reals{n\times d}$, 
  respectively, of variance $\sigma^{2}$, such that
  $\abs{\orig{a}} \leq 1$,  $\pnorm{2}{\orig{\bb}} \leq  1$,
  $\pnorm{2}{\orig{\xx }} \leq 1$, and $\pnorm{2}{\orig{\YY}} \leq 1$.
Let $\vv$ be an arbitrary vector in $\Reals{d}$.
If $a$, $\bb$, $\xx$, and $\YY$ are independent and
  $\sigma ^{2} \leq 1$, then 
\[
  \prob{}{\frac{\norm{\xx  + \YY \vv }_{\infty }
             }{
                \abs{a + \bb ^{T} \vv }} > x}
  \leq 
 \sqrt{\frac{2}{\pi }}
  \frac{1}{x}
\left(
\frac{ \sqrt{2} }{ \sigma }
+
  \sqrt{2 \ln max (n,2)} + \frac{1}{\sqrt{2\pi} \ln max (n,2)}
 \right),
\]
\end{lemma}
\begin{proof}
We begin by observing that $a + \bb ^{T} \vv $
  and each component of
  $\xx + \YY \vv$
  is a Gaussian random variable of variance
  $\sigma ^{2} (1 + \pnorm{2}{\vv}^{2})$ whose mean has
  absolute value at most $1 + \pnorm{2}{\vv }$,
and that all these variables are independent.
By Lemma~\ref{lem:maxGauss},
\[
  \expec{\xx,\YY}{\norm{\xx  + \YY \vv }_{\infty }}
\leq 
  1 + \pnorm{2}{\vv} + 
  \left(\sigma \sqrt{(1 + \pnorm{2}{\vv}^{2})} \right)
  \left(\sqrt{2 \ln max (n,2)} + \frac{1}{\sqrt{2\pi} \ln max (n,2)} \right).
\]
On the other hand,  Lemma~\ref{lem:distToPlane} implies
\begin{equation}\label{eqn:vecRatio}
  \prob{a,\bb}{\frac{1}{\abs{a + \bb ^{T} \vv}} > x}
\leq 
  \sqrt{\frac{2}{\pi }}
  \frac{1}{x \sigma \sqrt{1 + \pnorm{2}{\vv}^{2}}}.
\end{equation}
Thus, we can apply Corollary~\ref{lem:linear-expect} to show
\begin{align*}
 \prob{}{\frac{\norm{\xx  + \YY \vv }_{\infty }
             }{
                \abs{a + \bb ^{T} \vv }}
  > x}
& \leq 
 \sqrt{\frac{2}{\pi }}
 \frac{
  1 + \pnorm{2}{\vv} + 
  \left(\sigma \sqrt{1 + \pnorm{2}{\vv}^{2}} \right)
  \left(\sqrt{2 \ln max (n,2)} + \frac{1}{\sqrt{2\pi} \ln max (n,2)} \right)
   }{
     x \sigma \sqrt{1 + \pnorm{2}{\vv}^{2}}
   }\\
& =
 \sqrt{\frac{2}{\pi }}
  \frac{1}{x}
\left(
\frac{
  1 + \pnorm{2}{\vv} 
  }{
    \sigma \sqrt{1 + \pnorm{2}{\vv}^{2}}
  }
+
\frac{
  \left(\sigma  \sqrt{1 + \pnorm{2}{\vv}^{2}} \right)
  \left(\sqrt{2 \ln max (n,2)} + \frac{1}{\sqrt{2\pi} \ln max (n,2)} \right)
  }{
    \sigma \sqrt{1 + \pnorm{2}{\vv}^{2}}
  }
 \right)\\
& \leq 
 \sqrt{\frac{2}{\pi }}
  \frac{1}{x}
\left(
\frac{ \sqrt{2} }{ \sigma }
+
  \sqrt{2 \ln max (n,2)} + \frac{1}{\sqrt{2\pi} \ln max (n,2)}
 \right),
\end{align*}
where the last inequality follows from $(1+z)^{2}\leq 2 (1+z^{2})$,
  $\forall z \geq 0$.
\end{proof}

\section{Smoothed Analysis of Gaussian Elimination}\label{sec:smoothed}

We now combine the results from the previous sections to bound
  the smoothed precision needed 
  in the application of
  Gaussian elimination without pivoting
  to obtain
  solutions to linear systems accurate to $b$ bits.

\begin{theorem}[Smoothed precision of Gaussian elimination]\label{thm:gaussian}
For $n > e^{4}$,
  let $\orig{\AA}$ be an $n$-by-$n$ matrix for which $\pnorm{2}{\orig{\AA}} \leq 1$,
  and let $\AA$ be a Gaussian perturbation of $\orig{\AA}$ of variance 
  $\sigma ^{2} \leq  1/4$.
Then, the expected number of bits of precision necessary to 
  solve $\AA \xx = \bb$ to $b$
  bits of accuracy using Gaussian elimination without pivoting is
  at most
\[
 b + \frac{11}{2} \log_{2}n
   + 3 \log_{2} \left(\frac{1}{\sigma}\right)
   + \log_{2}(1+2\sqrt{n}\sigma) + \frac{1}{2}\log_{2}\log_{2} n
   +6.83
\]
\end{theorem}
\begin{proof}
By Wilkinson's theorem, we need the machine precision,
  $\epsilon _{mach}$, to satisfy
\begin{eqnarray*}
 5 \cdot  2^{b} n \rho _{\LL} (\AA) \rho _{\UU} (\AA) \kappa (\AA) \epsilon _{mach} 
  & \leq  & 1 \qquad \implies \\
2.33 + b + \log_{2} n + \log_{2} (\rho _{\LL} (\AA))  + 
  \max(0,\log_{2} (\rho _{\UU} (\AA))) +
  \log_{2} (\kappa (\AA))
& \leq  & \log_{2} (1/\epsilon _{mach} ).
\end{eqnarray*}
We will apply Lemma \ref{lem:linear-log} to bound these log terms.
Theorem~\ref{thm:U} tells us that
\[ \prob{}{\rho_U(\AA) > 1+x} \le 
   \frac{1}{\sqrt{2 \pi}}\frac{n (n+1)}{x \sigma} .
\]
To put this inequality into a form to which Lemma~\ref{lem:linear-log} may be applied,
  we set 
\[
y = x \left(1 + \frac{\sqrt{2 \pi} \sigma}{n (n+1)} \right),
\]
to obtain
\[
 \prob{}{\rho_U(\AA) > y} \le 
\left(\frac{1}{\sqrt{2\pi }}\frac{n (n+1)}{\sigma } + 1 \right)\frac{1}{y}.
\]
By Lemma \ref{lem:linear-log}, 
\begin{align*}
  \expec{}{ \max(0, \log_{2} \rho _{\UU} (\AA))} 
& \leq 
\log_{2} 
\left(\frac{1}{\sqrt{2\pi }}\frac{n (n+1)}{\sigma } + 1 \right) +
\log_{2}e\\
& \leq 
\log_{2} 
\left(n (n+1) + \sigma \sqrt{2 \pi } \right) +
\log_{2}\left(\frac{1}{\sigma}\right) +
\log_{2}\left(\frac{e}{\sqrt{2 \pi }} \right)\\
& \leq 
\log_{2} 
\left(1.02 n^{2} \right) +
\log_{2}\left(\frac{1}{\sigma}\right) +
\log_{2}\left(\frac{e}{\sqrt{2 \pi }} \right)\\
& \leq 
 2\log_{2}n + \log_{2}\left(\frac{1}{\sigma}\right) + 0.15
  ,
\end{align*}
where in the second-to-last inequality, we used the assumptions $n\geq e^{4}$ and
$\sigma \leq 1/2$.
In the last inequality, we numerically computed
  $\log_{2} (1.02 e / \sqrt{2 \pi}) < 0.15$.

Theorem~\ref{thm:L} and Lemma \ref{lem:linear-log} imply
\begin{align*}
  \expec{}{  \log_{2} \rho _{\LL} (\AA)}
& \leq 
\log_{2}\left(\sqrt{\frac{2}{\pi }}
  n^{2}
\left(
\frac{ \sqrt{2} }{ \sigma }
+
  \sqrt{2 \ln n} + \frac{1}{\sqrt{2\pi} \ln n}
 \right) \right) + \log_{2} e\\
  &\leq
 2 \log_{2}n
   + \log_{2}\left(\frac{1}{\sigma} + 
    \sqrt{\ln{n}}\left(1 + \frac{1}{2\sqrt{\pi } \ln n} \right)
  \right) +\log_{2} \left(\frac{2e}{\sqrt{\pi }} \right)\\
  & =
 2 \log_{2}n + \log_{2}\left(\frac{1}{\sigma}\right) 
  + \log_{2}\sqrt{\ln n}   + \log_{2}
   \left(\frac{1}{\sqrt{\ln  n}} + \sigma \left(1 + \frac{1}{2\sqrt{\pi } \ln n} \right)
  \right) +\log_{2} \left(\frac{2e}{\sqrt{\pi }} \right)\\
\intertext{using~$\sigma\leq\frac{1}{2}$ and $n > e^{4}$,}
& \leq 
 2 \log_{2}n + \log_{2}\left(\frac{1}{\sigma}\right) 
  + \frac{1}{2}\log_{2}\log_{2}n   + \log_{2}
   \left(1 + \frac{1}{16\sqrt{\pi }}
  \right) +\log_{2} \left(\frac{2e}{\sqrt{\pi }} \right)\\
  &\leq
 2 \log_{2}n
   + \log_{2}\left(\frac{1}{\sigma}\right)
   + \frac{1}{2}\log_{2}\log_{2} n + 1.67,
\end{align*}
as $\log_{2} (1+ 1/16\sqrt{\pi}) + \log_{2} (2e/\sqrt{\pi}) < 1.67$.
Theorem~\ref{thm:sigmaSmall} and Lemma \ref{lem:linear-log},
  along with the observation that $\log_{2} (2.35 e) < 2.68$,
 imply
\[
  \expec{}{  \log_{2} \pnorm{2}{\AA^{-1}}} \leq 
  \frac{1}{2}\log_{2}n + \log_{2}\left(\frac{1}{\sigma}\right) +  2.68.
\]
Finally,
\[
  \expec{}{  \log_{2} (\pnorm{2}{\AA})} \leq 
  \log_{2} (1+2 \sqrt{n} \sigma )
\]
follows from the well-known facts that the
  expectation of $\pnorm{2}{\AA - \orig{\AA}} $
  is at most $2 \sqrt{n} \sigma$
  (\textit{c.f.,} \cite{Seginer})
  and that
  $\expec{}{\log_{2} (X)} \leq \log_{2} \expec{}{X}$
  for every positive random variable $X$.
Thus,
  the expected number of digits of precision needed is at most
\[
 b + \frac{11}{2} \log_{2}n
   + 3 \log_{2} \left(\frac{1}{\sigma}\right)
   + \log_{2}(1+2\sqrt{n}\sigma) + \frac{1}{2}\log_{2}\log_{2} n
   +6.83.
\]
\end{proof}

The following conjecture would further improve the coefficient
  of $\log (1/ \sigma )$.

\begin{conjecture}\label{conj:combiningAll}
Let $\orig{\AA}$ be a $n$-by-$n$ matrix for which $\pnorm{2}{\orig{\AA}} \leq 1$,
  and let $\AA$ be Guassian perturbation of $\orig{\AA}$
  of variance $\sigma ^{2} \leq  1$.
Then 
\[
  \prob{}{\rho _{\LL} (\AA) \rho _{\UU} (\AA) \kappa (\AA) > x}
\leq \frac{ n^{c_{1}} \log^{c_{2}} (x)}{x \sigma},  
\]
for some constants $c_{1}$ and $c_{2}$.
\end{conjecture}

\section{Zero-preserving perturbations of symmetric matrices with diagonals}\label{sec:symmetric}

Many matrices that occur in practice are symmetric and 
  sparse.
Moreover, many matrix algorithms take advantage
  of this structure.
Thus, it is natural to study the smoothed analysis of algorithms
  under perturbations that respect symmetry and non-zero structure.
In this section, we study the condition numbers and growth factors
  of Gaussian elimination without pivoting of symmetric matrices
  under perturbations that only alter their diagonal and non-zero entries.

\begin{definition}[Zero-preserving perturbations]\label{def:zeroPert}
Let $\orig{\TT}$ be a matrix.
We define the \emph{zero-preserving perturbation of $\orig{\TT}$ 
  of variance $\sigma ^{2}$} to be the matrix $\TT$
  obtained by adding independent
  Gaussian random variables of mean 0
  and variance $\sigma ^{2}$ to the non-zero entries of $\orig{\TT}$.
\end{definition}
 
Throughout this section, when we express a symmetric
  matrix $\AA$ as $\TT + \DD + \TT^{T}$, we
  mean that $\TT$ is lower-triangular with zeros on the diagonal 
  and $\DD$ is a diagonal matrix.
By making a zero-preserving perturbation to $\orig{\TT}$, 
  we preserve the symmetry of the matrix.
The main results of this section are that the smoothed condition number
  and growth factors of symmetric matrices 
  under zero-preserving perturbations to $\TT$ and diagonal
  perturbations to $\DD$ 
  have distributions similar those proved
  in Sections~\ref{sec:cond} and~\ref{sec:growth} 
  for
  dense matrices under dense perturbations.

\subsection{Bounding the condition number}\label{ssec:symCond}

We begin by recalling that the singular values and vectors
  of symmetric matrices are the eigenvalues and eigenvectors.

\begin{lemma}\label{lem:symSigma1}
Let $\orig{\AA} = \orig{\TT} + \orig{\DD} + \orig{\TT}^{T}$ be an 
  arbitrary $n$-by-$n$ symmetric matrix.
Let $\TT$ be a zero-preserving perturbation of $\orig{\TT}$ of
  variance $\sigma ^{2}$, let $\GG_{D}$ be a diagonal matrix
  of independent Gaussian random variables of variance $\sigma ^{2}$
  and mean 0 that are independent of $\TT$, and let
  $\DD = \orig{\DD} + \GG_{D}$.
Then, for $\AA = \TT + \DD + \TT^{T}$,
\[
  \prob{}{\pnorm{2}{\AA^{-1}} \geq x }
  \leq 
  \sqrt{\frac2\pi}\frac{n^{3/2}}{x\sigma }.
\]
\end{lemma}

\begin{proof}
By Proposition \ref{pro:avemax},
\begin{eqnarray*}
  \prob{\TT,\GG_{D}}{\pnorm{2}{(\TT + \DD + \TT^{T})^{-1}} \geq x }
& \leq  &
 \max _{\TT}
  \prob{\GG_{D}}{\pnorm{2}{((\TT + \orig{\DD} + \TT^{T}) + \GG_{D})^{-1}} \geq x }.
\end{eqnarray*}
The proof now follows from Lemma~\ref{lem:diagSigma1}, taking
  $\TT + \orig{\DD} + \TT^{T}$ as the base matrix.
\end{proof}

\begin{lemma}\label{lem:diagSigma1}
Let $\orig{\AA}$ 
  be an arbitrary $n$-by-$n$ symmetric matrix,
  let $\GG_{D}$ be a diagonal matrix
  of independent Gaussian random variables of variance $\sigma ^{2}$
  and mean 0, and let
  $\AA = \orig{\AA} + \GG_{D}$.  
Then,  
\[
  \prob{}{\pnorm{2}{\AA^{-1}} \geq x }
  \leq 
  \sqrt{\frac{2}{\pi}}\frac{n^{3/2}}{x\sigma }.
\]
\end{lemma}

\begin{proof}
Let $\vs{x}{1}{n}$ be the diagonal entries of $\GG_{D}$, and
  let 
\begin{eqnarray*}
   g & = & \frac{1}{n} \sum _{i=1}^{n} x _{i}, \mbox{ and }\\
  y _{i} & = & x _{i} - g .
\end{eqnarray*}
Then,
\begin{eqnarray*}
  \prob{\vs{y}{1}{n},g}{\pnorm{2}{(\orig{\AA} + \GG_{D})^{-1}} \geq x }
& = &
  \prob{\vs{y}{1}{n},g}{\pnorm{2}{(\orig{\AA} + \textrm{diag} (\vs{y }{1}{n})
   + g \II)^{-1}} \geq x }\\
& \leq &
  \max _{\vs{y}{1}{n}}
  \prob{g}{\pnorm{2}{(\orig{\AA} + \textrm{diag} (\vs{y }{1}{n})
   + g \II)^{-1}} \geq x },
\end{eqnarray*}
where the last inequality follows from Proposition \ref{pro:avemax}.
The proof now 
 follows from Proposition~\ref{pro:coGaussians} and Lemma~\ref{lem:eyeSigma1}.
\end{proof}

\begin{proposition}\label{pro:coGaussians}
Let $\vs{X}{1}{n}$ be independent Gaussian random variables
  of variance $\sigma ^{2}$ with means $\vs{a}{1}{n}$, respectively.
Let 
\begin{eqnarray*}
   G & = & \frac{1}{n} \sum _{i=1}^{n} X _{i}, \mbox{ and }\\
   Y _{i} & = & X _{i} - G .
\end{eqnarray*}
Then, $G $ is a Gaussian random variable of variance
  $\sigma ^{2}/n$ with mean $(1/n) \sum a _{i}$,
  independent of $\vs{Y}{1}{n}$.
\end{proposition}

\begin{lemma}\label{lem:eyeSigma1}
Let $\orig{\AA}$ be an arbitrary $n$-by-$n$ symmetric matrix,
  and let $G$ be a
  Gaussian random variable of mean $0$ and variance $\sigma ^{2}/n$.
Let $\AA = \orig{\AA} + G\II$.  
Then,  
\[
  \prob{\AA}{\pnorm{2}{\AA^{-1}} \geq x }
  \leq 
  \sqrt{\frac{2}{\pi}}\frac{n^{3/2}}{ x \sigma}.
\]
\end{lemma}
\begin{proof}
Let $\vs{\lambda }{1}{n}$ be the eigenvalues of $\orig{\AA}$.
Then,
\[
  \pnorm{2}{(\orig{\AA} + G\II)^{-1}}^{-1}
 =
  \min _{i} \abs{\lambda _{i} + G}.
\]
Thus,
\begin{align*}
\prob{\AA}{\pnorm{2}{\AA^{-1}} \geq x }= \prob{G}{\min _{i} \abs{\lambda _{i} - G} < \frac{1}{x} }
   \leq  \sum_{i} \prob{G}{ \abs{\lambda _{i} - G} < \frac{1}{x} }
   \leq  \sum_{i} \sqrt{\frac{2}{\pi}}\frac{\sqrt{n} }{x \sigma }
   \leq  \sqrt{\frac{2}{\pi}}\frac{n^{3/2} }{x \sigma },
\end{align*}
where
  the second-to-last inequality follows from 
  Lemma~\ref{lem:distToPlane} for  $\Reals{1}$.
\end{proof}

As in Section~\ref{sec:cond}, we can now prove:
\begin{theorem}[Condition number of symmetric matrices]\label{thm:symCond}
Let $\orig{\AA} = \orig{\TT} + \orig{\DD} + \orig{\TT}^{T}$ be an 
  arbitrary $n$-by-$n$ symmetric matrix satisfying $\pnorm{2}{\orig{\AA}} \leq \sqrt{n}$.
Let $\sigma ^{2} \leq 1$, 
  let $\TT$ be a zero-preserving perturbation of $\orig{\TT}$ of
  variance $\sigma ^{2}$, let $\GG_{D}$ be a diagonal matrix
  of independent Gaussian random variables of variance $\sigma ^{2}$
  and mean $0$ that are independent of $\TT$, and let
  $\DD = \orig{\DD} + \GG_{D}$.
Then, for $\AA = \TT + \DD + \TT^{T}$,
\[
  \prob{}{\kappa (\AA) \geq x }
  \leq 
6\sqrt{\frac{2}{\pi}}\frac{n^{7/2}}{x\sigma}
  \left(1 + \sqrt{2\ln (x)/ 9n } \right)
\]
\end{theorem}
\begin{proof}
As in the proof of Theorem~\ref{thm:kappaSmall}, 
  we can apply the techniques used in the proof 
  of~\cite[Theorem II.7]{DavidsonSzarek},
  to show
\[
    \prob{}{\pnorm{2}{\orig{\AA} - \AA} \geq 2\sqrt{n} + k} < e^{-k^{2}/2}.
\]
The rest of the proof follows the 
  outline of the proof of Theorem~\ref{thm:kappaSmall},
  using Lemma~\ref{lem:symSigma1} instead of 
  Theorem~\ref{thm:sigmaSmall}.
\end{proof}

\subsection{Bounding entries in $\UU$}\label{ssec:symU}
In this section, we will prove:
\begin{theorem}[$\rho _{\UU} (\AA)$ of symmetric matrices]\label{thm:symU}
Let $\orig{\AA} = \orig{\TT} + \orig{\DD} + \orig{\TT}^{T}$ be an 
  arbitrary $n$-by-$n$ symmetric matrix satisfying $\pnorm{2}{\orig{\AA}} \leq 1$.
Let $\sigma ^{2} \leq 1$,
  let $\TT$ be a zero-preserving perturbation of $\orig{\TT}$ of
  variance $\sigma ^{2}$, let $\GG_{D}$ be a diagonal matrix
  of independent Gaussian random variables of variance $\sigma ^{2}$
  and mean $0$ that are independent of $\TT$, and let
  $\DD = \orig{\DD} + \GG_{D}$.
Then, for $\AA = \TT + \DD + \TT^{T}$,
\[
  \prob{}{\rho _{\UU} (\AA)
          > 1 + x}
   \leq \frac{2}{7}\sqrt{\frac2\pi}\frac{n^{3}}{x\sigma}
\]
\end{theorem}

\begin{proof}
We proceed as in the proof of Theorem~\ref{thm:U}.
For $k$ between $2$ and $n$, we define $\uu$,
  $\aa$, $\bb$ and $\CC $ as in the proof of Theorem
  \ref{thm:U}.
By \eqref{eqn:rhoU}
\begin{align*}
 \frac{\pnorm{1}{\uu}}{\pnorm{\infty}{\AA}}
   \leq 1+\pnorm{1}{\left(\CC^{T} \right)^{-1}\bb}
   \leq 1+\sqrt{k-1}\pnorm{2}{\bb^{T}\CC^{-1}} 
   \leq 1+\sqrt{k-1}\pnorm{2}{\bb}\pnorm{2}{\CC^{-1}}.
\end{align*}
Hence
\begin{align*}
 \prob{}{\frac{\pnorm{1}{\uu }}{\norm{\AA}_\infty} > 1+x}
   &\leq \prob{}
     {\pnorm{2}{\bb}\pnorm{2}{\CC^{-1}} >
     \frac{x}{\sqrt{k-1}}} \\
   &\leq \expec{}{\pnorm{2}{\bb}}
     \sqrt{\frac2\pi}\frac{(k-1)^{2}}{x\sigma} ,
  & \text{ by Lemmas~\ref{lem:symSigma1}
    and~\ref{lem:linear-expect},}
\\
   &\leq \sqrt{1+j\sigma^2}\sqrt{\frac2\pi} \frac{(k-1)^{2}}{x\sigma} ,
 & \text{where~$j$ is the number of non-zeros in~$\bb$,}\\
   &\leq \sqrt{\frac2\pi} \frac{\sqrt{k}(k-1)^{2}}{x\sigma}.
\end{align*}
Applying a union bound over~$k$,
\begin{align*}
 \prob{}{\rho_U(\AA) > x}
  \leq \sqrt{\frac2\pi}\frac{1}{x\sigma}
     \sum_{k=2}^n \sqrt{k}(k-1)^{2} \leq \frac{2}{7}\sqrt{\frac2\pi}\frac{n^{7/2}}{x\sigma}.
\end{align*}

\end{proof}

\subsection{Bounding entries in $\LL$}\label{ssec:symL}

As in Section~\ref{ssec:rhoL}, we derive a bound on
  the growth factor of $\LL$.
As before, we will show that it is unlikely that
  $\AA^{(k-1)}_{j,k}$ is large while $\AA^{(k-1)}_{k,k}$
  is small.
However, our techniques must differ from those used
  in Section~\ref{ssec:rhoL}, as
  the proof in that section made critical use of the
  independence of $\AA_{k,1:(k-1)}$ and $\AA_{1:(k-1),k}$.

\begin{theorem}[$\rho _{\LL} (\AA)$ of symmetric matrices]\label{thm:symL}
Let $\sigma ^{2} \leq 1$ and $n\geq 2$.
Let $\orig{\AA} = \orig{\TT} + \orig{\DD} + \orig{\TT}^{T}$ be an 
  arbitrary $n$-by-$n$ symmetric matrix satisfying 
  $\pnorm{2}{\orig{\AA}} \leq 1$.
Let $\TT$ be a zero-preserving perturbation of $\orig{\TT}$ of
  variance $\sigma ^{2}$, let $\GG_{D}$ be a diagonal matrix
  of independent Gaussian random variables of variance $\sigma ^{2} \leq 1$
  and mean $0$ that are independent of $\TT$, and let
  $\DD = \orig{\DD} + \GG_{D}$.
Let $\AA = \TT + \DD + \TT^{T}$.
Then, 
\[
\forall x \geq \sqrt{\frac{2}{\pi }}\frac{1}{\sigma^{2}},
\quad 
  \prob{}{\rho _{\LL} (\AA)
          > x}
  \leq 
\frac{3.2 n^{4}}{x \sigma^{2}}
    \ln^{3/2} \left( e \sqrt{\frac{\pi}{2}} x \sigma^{2}
   \right).
\]
\end{theorem}
\begin{proof}
Using Lemma~\ref{lem:bigLk}, 
  we obtain for all $k$
\begin{align*}
\prob{}{\exists j > k : \abs{\LL_{j,k}} > x}
& \leq 
\prob{}{\pnorm{2}{\LL_{(k+1):n,k}} > x}  \leq 
\frac{3.2 n^{2}}{x \sigma^{2}}
    \ln^{3/2} \left(
   e \sqrt{\frac{\pi}{2}} x \sigma^{2}
   \right).
\end{align*}
Applying a union bound over the choices for $k$,
  we then have
\[
  \prob{}{\exists j, k : \abs{\LL_{j,k}} > x}
  \leq 
\frac{3.2 n^{3}}{x \sigma^{2}}
    \ln^{3/2} \left(
   e \sqrt{\frac{\pi}{2}} x \sigma^{2}
   \right).
\]
The result now follows from the fact that
 $\norm{\LL}_{\infty }$ is at most $n$ times the
  largest entry in $\LL$.
\end{proof}

\begin{lemma}\label{lem:bigLk}
Under the conditions of Theorem~\ref{thm:symL},
\[
\forall x \geq \sqrt{\frac{2}{\pi }}\frac{1}{\sigma^{2}},
\quad 
\prob{}{\pnorm{2}{\LL_{(k+1):n,k}} > x}
  \leq 
\frac{3.2 n^{2}}{x \sigma^{2}}
    \ln^{3/2} \left(
   e \sqrt{\frac{\pi}{2}} x \sigma^{2}
   \right).
\]
\end{lemma}
\begin{proof}
We recall that
\[
  \LL_{k+1:n,k} = \frac{\AA_{k+1:n,k} - \AA_{k+1:n,1:k-1}\AA_{1:k-1,1:k-1}^{-1}\AA_{1:k-1,k}}
    {\AA_{k,k} - \AA_{k,1:k-1}\AA_{1:k-1,1:k-1}^{-1}\AA_{1:k-1,k}}
\]
Because of the symmetry of~$\AA$, 
  $\AA_{k,1:k-1}$ is the same as~$\AA_{1:k-1,k}$, so we
  can no longer use the proof technique that worked in Section~\ref{ssec:rhoL}.
Instead, we will bound the tails of the numerator and denominator separately,
  exploiting the fact that only the denominator depends upon $\AA_{k,k}$.

Consider the numerator first.
Setting $\vv = \AA_{1:k-1,1:k-1}^{-1}\AA_{1:k-1,k}$, the numerator
  can be written
  $\AA_{k+1:n, 1:k} \left({-\vv} \atop { 1} \right)$.
We will now prove that for all $x \geq 1/\sigma$,
\begin{equation}\label{eqn:bigLk}
  \prob{\substack{\AA_{k+1:n, 1:k}\\
     A_{1:k-1, 1:k}}}
  {\norm{\AA_{k+1:n, 1:k} \left({-\vv} \atop { 1} \right)}_{\infty } > x}
  \leq 
\sqrt{\frac{2}{\pi }}
\left(
  \frac{2 n^{2} (1 + \sigma \sqrt{2 \ln (x\sigma )}) + n}
       {x \sigma }
  \right).
\end{equation}

Let 
\begin{equation}
  c = \frac{1}{1 + \sigma \sqrt{2 \ln (x\sigma )}},
\end{equation}
which implies $\frac{1-c}{c\sigma } = \sqrt{2 \ln (x\sigma )}$.
It suffices to prove \eqref{eqn:bigLk} for all $x$ for which the right-hand side
  is less than $1$. 
Given that $x \geq 1/\sigma$,
  it suffices to consider $x$ for which $cx \geq 2$
  and $x \sigma \geq 2$.

We use the parameter $c$ to divide the probability as follows:
\begin{align}
&   \prob{\substack{\AA_{k+1:n, 1:k}\\
     A_{1:k-1, 1:k}}}
  {\pnorm{\infty}{\AA_{k+1:n, 1:k} \left({-\vv} \atop { 1} \right)} > x} \nonumber\\
& \quad \leq 
  \prob{\AA_{1:(k-1), 1:k}}
  {\pnorm{2}{\left({-\vv} \atop { 1} \right)} > c x} \label{eqn:bigLk1}\\
& \qquad +
  \prob{\AA_{k+1:n, 1:k}}
  {\pnorm{\infty}{\AA_{k+1:n, 1:k} \left({-\vv \atop { 1} }\right)} > 
  \frac{1}{c} \pnorm{2}{\left(-\vv \atop { 1} \right)}
  \Bigg|
 \pnorm{2}{\left(-\vv \atop { 1} \right)} \leq  cx
 } \label{eqn:bigLk2}
\end{align}

To evaluate \eqref{eqn:bigLk2}, we note that
  once $\vv$ is fixed, each component of 
  $\AA_{k+1:n, 1:k} \left({-\vv} \atop { 1} \right)$
  is a Gaussian random variable of variance
$\pnorm{2}{\left({-\vv} \atop { 1} \right)}^{2}\sigma^{2}$
  and mean at most
  $\pnorm{2}{\orig{\AA}_{k+1:n, 1:k} \left({-\vv} \atop { 1} \right)}
  \leq \pnorm{2}{\left({-\vv} \atop { 1} \right)}$.
So,
\[
\pnorm{\infty}{\AA_{k+1:n, 1:k} \left({-\vv} \atop { 1} \right)} > 
  \frac{1}{c} \pnorm{2}{\left({-\vv} \atop { 1} \right)}
\]
 implies one of the Gaussian random variables differs from its mean
  by more than $(1/c -1)/\sigma $ times it standard deviation, and we can therefore
  apply Lemma~\ref{lem:gaussTail} and a union bound to derive
\[
 \eqref{eqn:bigLk2} 
\leq 
 \sqrt{\frac{2}{\pi }}
  \frac{n e^{-\frac{1}{2} \left(\frac{1-c}{c\sigma } \right)^{2}}}
   {\frac{1-c}{c \sigma }  }
=
\sqrt{\frac{2}{\pi }}
\frac{n}{x \sigma \sqrt{2 \ln (x \sigma )}}
.
\]
To bound \eqref{eqn:bigLk1}, we note that
 Lemma~\ref{lem:symSigma1} and Corollary~\ref{cor:chiComb} 
  imply
\[
  \prob{\AA_{1:(k-1), 1:k}}
  {\pnorm{2}{\AA_{1:k-1,1:k-1}^{-1}\AA_{1:k-1,k}} > y}
 \leq 
 \sqrt{\frac{2}{\pi }}
  \frac{n^{2}}{y \sigma }
,
\]
and so
\begin{align*}
  \prob{\AA_{1:(k-1), 1:k}}
  {\pnorm{2}{\left({-\vv} \atop { 1} \right)} > c x}
& \leq 
  \prob{\AA_{1:(k-1), 1:k}}
  {\pnorm{2}{\AA_{1:k-1,1:k-1}^{-1}\AA_{1:k-1,k}} > cx - 1}\\
& \leq 
 \sqrt{\frac{2}{\pi }}
  \frac{n^{2}}{(cx - 1) \sigma }\\
& =
 \sqrt{\frac{2}{\pi }}
  \frac{n^{2}}{(cx\sigma (1-1/cx))}\\
& =
 \sqrt{\frac{2}{\pi }}
  \frac{n^{2} (1 + \sigma \sqrt{2 \ln (x\sigma )})}
       {x \sigma \left(1-1/cx\right)} \\
& \leq 
 \sqrt{\frac{2}{\pi }}
  \frac{2 n^{2} (1 + \sigma \sqrt{2 \ln (x\sigma )})}
       {x \sigma }, \text{ by $cx\geq 2$.}
\end{align*}
So,
\begin{align}
  \prob{\substack{\AA_{k+1:n, 1:k}\\
     A_{1:k-1, 1:k}}}
  {\norm{\AA_{k+1:n, 1:k} \left({-\vv} \atop { 1} \right)}_{\infty } > x}
& \leq 
\sqrt{\frac{2}{\pi }}
\left(
\frac{n}{x \sigma \sqrt{2 \ln (x \sigma )}}
+
  \frac{2 n^{2} (1 + \sigma \sqrt{2 \ln (x\sigma )})}
       {x \sigma }
 \right) \notag \\
& \leq 
\sqrt{\frac{2}{\pi }}
\left(
  \frac{2 n^{2} \left(1 + \sigma \sqrt{2 \ln (x\sigma )} \right) + n}
       {x \sigma }
 \right), \label{eqn:Lnumerator}
\end{align}
by the assumption $x \sigma \geq 2$, which proves \eqref{eqn:bigLk}.

As for the denominator, 
  we note that~$\AA_{k,k}$ is independent of 
  all other terms, and hence
\begin{equation}\label{eqn:Ldenominator} 
\prob{}{\abs{\AA_{k,k} - \AA_{k,1:k-1}\AA_{1:k-1,1:k-1}^{-1}\AA_{1:k-1,k}} < 1/x}
       \leq \sqrt{\frac{2}{\pi}}\frac{1}{x\sigma}, 
\end{equation}
by Lemma~\ref{lem:distToPlane}.
Applying Corollary~\ref{cor:loglinComb} 
  with 
\[
 \alpha = \sqrt{\frac{2}{\pi}} \left( 2n^{2} + n \right)
\qquad 
  \beta = \frac{4 n^{2} \sigma }{\sqrt{\pi}} 
\qquad 
  \gamma = \sqrt{\frac{2}{\pi }}
\]
  to combine \eqref{eqn:Lnumerator} with \eqref{eqn:Ldenominator},
  we derive the bound 
\begin{multline*}
\frac{2}{\pi x \sigma^{2}}
\left(2n^{2} + n
+
\left(
 \left(2 + 4 \sqrt{2} \sigma /3 \right)n^{2}
 + n
 \right)
    \ln^{3/2} \left(
   \sqrt{\pi /2} x \sigma^{2}
   \right)
 \right)\\
\leq 
\frac{2 n^{2}}{\pi x \sigma^{2}}
 \left(3 + 4 \sqrt{2} \sigma /3 \right)
   \left( \ln^{3/2} \left(
   \sqrt{\pi /2} x \sigma^{2}
   \right) + 1 \right)\\
\leq 
\frac{3.2 n^{2}}{x \sigma^{2}}
    \ln^{3/2} \left(
   e \sqrt{\pi /2} x \sigma^{2}
   \right),
\end{multline*}
as $\sigma \leq 1$.
\end{proof}

\section{Conclusions and open problems}\label{sec:conclusions}

\subsection{Generality of results}
In this paper, we have presented bounds on the smoothed values of the 
   condition number and growth factors assuming the input matrix is subjected
   to a slight Gaussian perturbation.
We would like to point out here that 
   our results can be extended to some other families of perturbations.

With the exception of the proof of Theorem~\ref{thm:sigmaSmall},
  the only properties of Gaussian random
  vectors that we used in Sections~\ref{sec:cond} 
  and~\ref{sec:growth} are
\begin{enumerate}
\item [1.] there is a constant $c$ for which the probability
  that a Gaussian random vector has distance less than $\epsilon $
  to a hyperplane is at most $c \epsilon $, and
\item [2.] it is exponentially unlikely that a Gaussian random vector
 lies far from its mean.
\end{enumerate}
Moreover, a result similar to Theorem~\ref{thm:sigmaSmall} but with
  an extra factor of $d$ could be proved using just
  fact 1.

In fact, results of a character similar to ours would still hold if
  the second condition were reduced to a polynomial probability.
Many other families of perturbations share these properties.
For example, similar results would hold if we let
  $\AA = \orig{\AA} +\UU $, where $\UU $ is a matrix of variables independently
  uniformly chosen in $[-\sigma ,\sigma ]$, or if
  $\AA = \orig{\AA} + \SS$, where the columns of $\SS$ are chosen
  uniformly among those vectors of norm at most $\sigma $.

\subsection{Counter-Examples}
The results of sections~\ref{sec:cond} and~\ref{sec:growth} do not
  extend to zero-preserving perturbations for non-symmetric matrices.
For example, the following matrix remains ill-conditioned under
  zero-preserving perturbations.
\[
\begin{array}{ccccc}
     1&    -2&     0&     0&     0\\
     0&     1&    -2&     0&     0\\
     0&     0&     1&    -2&     0\\
     0&     0&     0&     1&    -2\\
     0&     0&     0&     0&     1
\end{array}
\]
A symmetric matrix that remains ill-conditioned under zero-preserving
  perturbations that do not alter the diagonal can be obtained by
  locating the above matrix in the upper-right quadrant, and its
  transpose in the lower-left quadrant:
\[
\begin{array}{cccccccccc}
     0&     0&     0&     0&     0&     1&    -2&     0&     0&     0\\
     0&     0&     0&     0&     0&     0&     1&    -2&     0&     0\\
     0&     0&     0&     0&     0&     0&     0&     1&    -2&     0\\
     0&     0&     0&     0&     0&     0&     0&     0&     1&    -2\\
     0&     0&     0&     0&     0&     0&     0&     0&     0&     1\\
     1&     0&     0&     0&     0&     0&     0&     0&     0&     0\\
    -2&     1&     0&     0&     0&     0&     0&     0&     0&     0\\
     0&    -2&     1&     0&     0&     0&     0&     0&     0&     0\\
     0&     0&    -2&     1&     0&     0&     0&     0&     0&     0\\
     0&     0&     0&    -2&     1&     0&     0&     0&     0&     0\\
\end{array}
\]

The following matrix maintains large
 growth factor under zero-preserving perturbations,
  regardless of whether partial pivoting or no pivoting is used.
\[
\begin{array}{cccccc}
    1.1  &       0  &       0  &       0  &       0  &  1\\
   -1  &  1.1  &       0  &       0  &       0  &  1\\
   -1  & -1  &  1.1  &       0  &       0  &  1\\
   -1  & -1  & -1  &  1.1  &       0  &  1\\
   -1  & -1  & -1  & -1  &  1.1  &  1\\
   -1  & -1  & -1  & -1  & -1  &  1\\
\end{array}
\]

These examples can be easily normalized to so that their 2-norms are
  equal to 1.

\subsection{Open Problems}

Questions that naturally follow from this work are:

\begin{itemize}
\item What is the probability that the perturbation of
  an arbitrary matrix has large growth factors under
  Gaussian elimination with partial pivoting?

\item What is the probability that the perturbation of
  an arbitrary matrix has large growth factors under
  Gaussian elimination with complete pivoting?

\item Can zero-preserving perturbations of symmetric matrices 
  have large growth factors under partial pivoting or under complete pivoting?

\item Can  zero-preserving perturbations of arbitrary matrices 
  have large growth factors under complete pivoting?
\end{itemize}
For the first question, we point out that experimental
  data of Trefethen and Bau \cite[p. 168]{TrefethenBau}
  suggest that the probability that the perturbation of an arbitrary
  matrix has large growth factor
  under partial pivoting may be exponentially smaller than
  without pivoting.
This leads us to conjecture:
\begin{conjecture}
Let $\orig{\AA}$ be an $n$-by-$n$ matrix for which
  $\pnorm{2}{\orig{\AA}} \leq 1$, and let $\AA$ be a Gaussian perturbation
  of $\orig{\AA}$ of variance $\sigma ^{2} \leq  1$.
Let $\UU$ be the upper-triangular matrix obtained from the LU-factorization
 of $\AA$ with partial pivoting.
There exist absolute constants $k_{1}$, $k_{2}$ and $\alpha $ for which
\[
  \prob{}{\maxnorm{\UU} / \maxnorm{\AA}
          > x + 1}
  \leq 
  n^{k_{1}} e^{-\alpha x^{k_{2}} \sigma }
\]
\end{conjecture}

Finally, we ask whether similar analyses can be performed for
  other algorithms of Numerical Analysis.
One might start by extending Smale's program by analyzing the
  smoothed values of other condition numbers.

\subsection{Recent Progress}\label{sub:recent}

Since the announcement of our result, Wschebor \cite{Wschebor}
  improved the smoothed bound on the condition number.

\begin{theorem}[Wschebor]\label{thm:Wschebor}
Let $\orig{\AA}$ be an $n \times n$~matrix 
  and let $\AA$ be a Gaussian perturbation of $\orig{\AA}$ of variance 
  $\sigma ^{2} \leq  1$.
Then,
\[
  \prob{}{\kappa (\AA) \geq x }
  \leq 
\frac{n}{x}\left(\frac{1}{4\sqrt{2\pi n}} +
7\left(5+\frac{4\pnorm{2}{\orig{\AA }}^{2} (1+\log n)}{\sigma^{2}n} \right)^{1/2} \right)
\]
\end{theorem}

When $\pnorm{2}{\orig{\AA }}\leq \sqrt{n}$, his result implies

\[
  \prob{}{\kappa (\AA) \geq x }
  \leq O\left(\frac{n\log n}{x\sigma } \right).
\]

We conjecture

\begin{conjecture}\label{conj:cond}
Let $\orig{\AA}$ be an $n \times n$~matrix satisfying $\pnorm{2}{\orig{\AA}}\leq \sqrt{n}$,
  and let $\AA$ be a Gaussian perturbation of $\orig{\AA}$ of variance 
  $\sigma ^{2} \leq  1$.
Then,
\[
  \prob{}{\kappa (\AA) 
          \geq 
          x
         }
  \leq 
   O \left(\frac{n}{x\sigma }\right).
\]
\end{conjecture}

\section{Acknowledgments}\label{sec:ack}
We thank Alan Edelman for
  suggesting the name ``smoothed analysis'', for
  suggesting we examine growth factors,
  and for his continuing support of our efforts.
We thank Juan Cuesta and Mario Wschebor for pointing out
  some mistakes in an early draft of this paper.
We thank Felipe Cucker for bringing Wschebor's paper \cite{Wschebor} to our attention.
Finally, we thank the referees for their extraordinary efforts and many
  helpful suggestions.

\bibliographystyle{alpha}
\bibliography{numer}

\appendix

\section{Gaussian random variables}\label{sec:gaussian}


\begin{lemma}\label{lem:gaussTail}
Let $X$ be a univariate Gaussian random variable with mean $0$ and
standard deviation $1$. Then for all $k\geq 1$,
\[
  \prob{}{X \geq k} 
    \leq \frac1{\sqrt{2\pi} } \frac{e^{-\frac12 k^2}}k.
\]
\end{lemma}
\begin{proof}
We have
\begin{align*}
  \prob{}{X \geq k}
    &= \frac1{\sqrt{2\pi}} \int_k^\infty e^{-\frac12 x^2}\,dx \\
    \intertext{putting~$t=\frac12 x^2$,}
    &= \frac{1}{\sqrt{2\pi}} \int_{\frac{1}{2} k^2}^\infty \frac{e^{-t}}{\sqrt{2t}} \,dt \\
    &\le \frac1{\sqrt{2\pi}} \int_{\frac12 k^2}^\infty \frac{e^{-t}}k \,dt \\
    &= \frac1{\sqrt{2\pi}} \frac{e^{-\frac12 k^2}}k .
\end{align*}
\end{proof}

\begin{lemma}\label{lem:distToPlane}
Let $\xx$ be a $d$-dimensional Gaussian random vector of
  variance $\sigma^{2}$, let $\tt$ be a unit vector, and let
  $\lambda$ be a real.
Then,
\[
  \prob{}{\abs{\form{\tt}{\xx} - \lambda } \leq \epsilon} \leq 
     \sqrt{\frac{2}{\pi }}\frac{\epsilon}{\sigma }.
\]
\end{lemma}

\begin{lemma}\label{lem:maxGauss}
Let $\vs{g}{1}{n}$ be 
  Gaussian random variables of mean 0
  and variance $1$.
Then,
\[
  \expec{}{\max_{i} \abs{g_{i}}} \leq 
  \sqrt{2\ln(\max (n,2))} + \frac{1}{\sqrt{2\pi}\ln(\max (n,2))}.
\]
\end{lemma}
\begin{proof}
For any $a\geq 1$,
\begin{align*}
  \expec{}{\max_{i} \abs{g_{i}}}
& = 
\int_{t = 0}^{\infty}
  \prob{}{\max_{i} \abs{g_{i}} \geq t} \diff{t}\\
& \leq 
\int_{t=0}^{a}
 1 \diff{t}
+
  \int_{a}^{\infty }
  n \prob{}{\abs{g_{1}} \geq t} \diff{t}\\
& \leq 
  a
 +
  \int_{a}^{\infty }
  n \frac{2}{\sqrt{2\pi}} \frac{e^{-\frac{1}{2}t^{2}}}{t} \diff{t}
&(\text{applying Lemma~\ref{lem:gaussTail},})
\\
& =
  a
 +
  \frac{2n}{\sqrt{2\pi}} \int_{a}^{\infty }
  \frac{e^{-\frac{1}{2}t^{2}}}{t^{2}} \diff{\left(\frac{1}{2}t^{2} \right)}
\\
& \leq 
  a
 +
  \frac{2n}{\sqrt{2\pi}}\frac{1}{a^{2}} \int_{a}^{\infty }
   e^{-\frac{1}{2}t^{2}}  \diff{\left(\frac{1}{2}t^{2} \right)}
\\
& =
  a
 +
  \frac{2n}{\sqrt{2\pi}}\frac{1}{a^{2}} e^{-\frac{1}{2}a^{2}}.
\end{align*}
Setting $a = \sqrt{2\ln(\max (n,2))}$, which is greater
  than 1 for all $n\geq 1$,
 we obtain the following upper bound on the expectation:
\[
\sqrt{2\ln(\max (n,2))} + \frac{2n}{\sqrt{2\pi }}\frac{1}{2\ln(\max
(n,2))}  \frac{1}{\max (n,2)}\leq 
\sqrt{2\ln(\max (n,2))} + \frac{1}{\sqrt{2\pi}\ln(\max (n,2))}.
\]

\end{proof}

\begin{lemma}[Expectation of reciprocal of the 1-norm of a Gaussian
vector]
\label{lem:expected1normGaussian}
Let $\orig{\aa}$ be an arbitrary column vector in $\Reals{n}$ for $n\geq 2$.
Let $\aa$ be a Gaussian perturbation of $\orig{\aa}$ of variance $\sigma^{2}$.
Then 
\[
\expec{}{\frac{1}{\norm{\aa}_{1}}} \leq 
  \frac{2}{n\sigma }
\]
\end{lemma}
\begin{proof} Let $\aa= (\vs{a}{1}{n})$. 
It is clear that the expectation of $1/\norm{\aa}_1$ is maximized if 
  $\orig{\aa} = \origin $, so we will make this assumption.
Without loss of generality,
  we also assume $\sigma^{2}=1$. 
For general $\sigma $, we can simply scale the bound by the factor
  $1/\sigma $.

Recall that the Laplace transform of a positive random variable X is defined by
\[
  \laplace{X}(t) = \expec{X}{e^{-tX}}
\]
and the expectation of the reciprocal of a random variable is simply the
  integral of its Laplace transform.

Let~$X$ be the absolute value of a standard normal random variable.
The Laplace transform of~$X$ is given by
\begin{eqnarray*}
  \laplace{X}(t)
    &=& \sqrt{\frac{2}{\pi}} \int_0^\infty e^{-tx} e^{-\frac12 x^2} \diff{x} \\
    &=& \sqrt{\frac{2}{\pi}}e^{\frac12 t^2}
      \int_0^\infty e^{-\frac12 (x+t)^2}\diff{x} \\
    &=& \sqrt{\frac{2}{\pi}} e^{\frac12 t^2}
      \int_t^\infty e^{-\frac12 x^2} \diff{x} \\
    &=& e^{\frac12 t^2}\erfc\left(\frac{t}{\sqrt{2}}\right).
\end{eqnarray*}
Taking second derivatives, and applying the inequality (\textit{c.f.} \cite[26.2.13]{AbramowitzStegun})
\[
  \frac{1}{\sqrt{2\pi }} \int_t^\infty e^{-\frac12 x^2} \diff{x}
\geq 
  \frac{e^{-\frac12 x^2}}{\sqrt{2\pi }} \frac{1}{x + 1/x},
\]
we find that
 $e^{\frac12 t^2}\erfc\left(\frac{t}{\sqrt{2}}\right)$ is convex.

We now set a constant $c = 2.4$ and set $\alpha$ to satisfy
\[
  1 - \frac{\sqrt{c/\pi}}{\alpha} =  
  e^{\frac12 (c/\pi )}\erfc\left(\frac{\sqrt{c/\pi }}{\sqrt{2}}\right).
\]
Numerically, we find that $\alpha \approx 1.9857 < 2$.

As $e^{\frac12 t^2}\erfc\left(\frac{t}{\sqrt{2}}\right)$ is convex,
  we have the upper bound
\[
 e^{\frac12 t^2}\erfc\left(\frac{t}{\sqrt{2}}\right)
\leq 
 1 - \frac{t}{\alpha}, \text{ for $0 \leq t \leq \sqrt{c/\pi }$.}
\]
For $t > \sqrt{c/\pi }$, we apply the upper bound
\[
 e^{\frac12 t^2}\erfc\left(\frac{t}{\sqrt{2}}\right)
\leq 
 \sqrt{\frac{2}{\pi}}\frac{1}{t},
\]
which follows from Lemma~\ref{lem:gaussTail}.

We now have
\begin{align*}
  \expec{}{\frac{1}{\norm{\aa}_1}}
    &= \int_0^\infty
      \left( e^{\frac{1}{2}t^2}\erfc(t/\sqrt{2}) \right)^n \diff{t} \\
    &\leq \int_0^{\sqrt{c/\pi}} \left( 1-\frac{t}{\alpha} \right)^n \diff{t}
      + \int_{\sqrt{c/\pi}}^\infty
          \left( \sqrt{\frac{2}{\pi}}\frac{1}{t} \right)^n \diff{t} \\
    &\leq \frac{\alpha}{n+1}+\sqrt{\frac{2}{\pi}}\frac{(2/c)^{(n-1)/2}}{n-1} \\
    &< \frac{2}{n+1}+\sqrt{\frac{2}{\pi}}\frac{(2/c)^{(n-1)/2}}{n-1} \\
    &\leq \frac{2}{n-1},
\end{align*}
for $n \geq 2$.
To verify this last equality, one can multiply through by $(n+1) (n-1)$ to obtain
\[
\sqrt{\frac{2}{\pi}} (n+1) (2/c)^{(n-1)/2} \leq 4,
\]
which one can verify by taking the derivitive of the left-hand side to 
  find the point where it is maximized, $n = (2 + \ln (5/6)) / \ln (6/5)$.
\end{proof}

\section{Random point on a sphere}\label{sec:sphere}
\begin{lemma}\label{lem:randSphere}
Let $d \geq 2$ and 
  let $\left(\vs{u}{1}{d} \right)$ be a 
  unit vector chosen uniformly at random in $\Reals{d}$.
Then,  for $c \leq 1$,
\[
  \prob{}{\abs{u_{1}} \geq \sqrt{\frac{c}{d}}}
\geq 
  \prob{}{\abs{G} \geq \sqrt{c}},
\]
where $G$ is a Gaussian random variable of variance $1$ and mean $0$.
\end{lemma}
\begin{proof}
We may obtain a random unit vector by choosing
  $d$ independent Gaussian random variables of variance $1$
  and mean $0$, $\vs{x}{1}{d}$, and setting
\[
  u_{i} = \frac{x_{i}}{\sqrt{x_{1}^{2} + \dotsb + x_{d}^{2}}}.
\]
We have
\begin{align*}
  \prob{}{u_{1}^{2} \geq \frac{c}{d}}
& = 
  \prob{}{\frac{x_{1}^{2}}{x_{1}^{2} + \dotsb + x_{d}^{2}} \geq \frac{c}{d}}\\
& = 
  \prob{}{\frac{(d-1) x_{1}^{2}}{x_{2}^{2}+ \dotsb + x_{d}^{2}} \geq
    \frac{(d-1)c}{d-c}}\\
& \geq 
  \prob{}{\frac{(d-1) x_{1}^{2}}{x_{2}^{2} + \dotsb + x_{d}^{2}} \geq c},
  \text{ since~$c\leq1$.}
\end{align*}
We now note that 
\[
   t_{d} 
\defeq
    \frac{\sqrt{(d-1)} x_{1}}
         {\sqrt{x_{2}^{2} + \dotsb + x_{d}^{2}}}
\]
is a random variable distributed according to
  the $t$-distribution with $d-1$ degrees of freedom.
The lemma now follows from the fact
(\textit{c.f.} \cite[Chapter 28, Section 2]{JohnsonKotzBala2} or
 \cite[26.7.5]{AbramowitzStegun}) that, for $c > 0$,
\[
  \prob{}{t_{d} > \sqrt{c}} \geq \prob{}{G > \sqrt{c}},
\]
and that the distributions of $t_{d}$ and $G$ are symmetric about the origin.
\end{proof}

\section{Combination Lemmas}\label{sec:combine}
\begin{lemma}\label{lem:combLemma}
Let~$A$ and~$B$ be two positive random variables.
Assume
  \begin{enumerate}
    \item $ \prob{}{A \ge x} \le f(x) $.
    \item $ \prob{}{B \ge x|A} \le g(x) $.
  \end{enumerate}
  where $g$~is monotonically decreasing and~$\lim_{x\to\infty}g(x)=0$.
Then,
  \[ \prob{}{AB \ge x} \le \int_0^\infty
      f\left(\frac xt\right) (-g'(t))\,dt
  \]
\end{lemma}
\begin{proof}
Let~$\mu_A$ denote the probability measure associated with~$A$.
We have
  
\begin{align*}
    \prob{}{AB \ge x}
      &= \int_0^\infty \prob{B}{B \ge x/s| A}\,d\mu_A(s) \\
      &\le \int_0^\infty g\left(\frac xs\right)\,d\mu_A(s), \\
    \intertext{integrating by parts,}
      &= \int_0^\infty \prob{}{A \ge s}
           \frac{d}{ds}g\left(\frac xs\right)\,ds \\
      &\le \int_0^\infty f(s) \frac{d}{ds}g\left(\frac xs\right)\,ds , \\
   \intertext{setting $ t = x/s$}
      &= \int_0^\infty f\left(\frac xt\right) (-g'(t))\,dt .
  \end{align*}
\end{proof}

\begin{corollary}[linear-linear]\label{cor:linComb}
Let~$A$ and~$B$ be two positive random variables.
Assume
  \begin{enumerate}
    \item $ \prob{}{A \ge x} \le \frac\alpha x $ and
    \item $ \prob{}{B \ge x| A} \le \frac\beta x $
  \end{enumerate}
  for some~$\alpha,\beta > 0$.
Then,
  \[ \prob{}{AB \ge x} \le \frac{\alpha\beta}x
       \left(1+\max\left(0,\ln\left(\frac x{\alpha\beta}\right)\right)\right) \]
\end{corollary}
\begin{proof}
As the probability of an event can be at most $1$,
\begin{align*}
   \prob{}{A \geq x} 
 &\leq 
  \min \left(\frac{\alpha }{x}, 1 \right)
 \defeq f (x), \mbox{ and}\\
   \prob{}{B \geq x} 
 &\leq 
  \min \left(\frac{\beta  }{x}, 1 \right)
  \defeq g (x).
\end{align*}
Applying Lemma~\ref{lem:combLemma} while observing
\begin{itemize}
\item $g' (t) = 0$ for $t \in \left[0,\beta  \right]$, and
\item $f (x/t) = 1$ for $t \geq x/\alpha $,
\end{itemize}
we obtain 
\begin{align*}
  \prob{}{A B \geq x}
& \leq 
 \int _{0}^{\beta } \frac{\alpha t}{x} \cdot 0 \diff{t}
+
\max \left(0, \int _{\beta }^{x/\alpha } 
  \frac{\alpha t}{x} \frac{\beta }{t^{2}} \diff{t}\right) 
+
 \int _{x/\alpha }^{\infty } \frac{\beta }{t^{2}} \diff{t}\\
& =
\max \left(0, \frac{\alpha \beta }{x} \int _{\beta }^{x/\alpha } \frac{d t}{t} \right) 
+
  \frac{\alpha \beta }{x}\\
& = \frac{\alpha \beta }{x}
\left(1 + \max \left(0, \ln \left(\frac{x}{\alpha \beta } \right) \right)\right),
\end{align*}
where the $\max$ appears in case $x/\alpha < \beta$.
\end{proof}

\begin{corollary}\label{cor:loglinComb}
Let~$A$ and~$B$ be two positive random variables.
If
  \begin{enumerate}
    \item $\forall x \geq 1/\sigma$,  $ \prob{}{A \ge x} \le
   \min \left(1,  \frac{\alpha+\beta\sqrt{\ln x\sigma}}
      {\sigma x} \right) $ and
    \item $ \prob{}{B \ge x| A} \le \frac\gamma {x\sigma} $
  \end{enumerate}
  for some~$\alpha \geq 1$ and  $\beta,\gamma,\sigma > 0$,
 then,
  \[
\forall x \geq \gamma /\sigma^{2}, \,
 \prob{}{AB \ge x}
    \leq
 \frac{\alpha\gamma}{x\sigma^2}
      \left(1 
        + \left(\frac{2\beta}{3\alpha} + 1 \right)
	  ln^{3/2}\left(\frac{x\sigma^2}{\gamma} \right)
          \right).
  \]
\end{corollary}
\begin{proof}
Define~$f$ and~$g$ by
\begin{align*}
  f(x) &\defeq
    \begin{cases}
      1 & \mbox{for~$x \leq \frac{\alpha}{\sigma}$}\\
      \frac{\alpha+\beta\sqrt{\ln x\sigma}}{x\sigma}
        & \mbox{for~$x > \frac{\alpha}{\sigma}$}
    \end{cases}\\
  g(x) &\defeq
    \begin{cases}
      1 & \mbox{for~$x \leq \frac{\gamma}{\sigma}$}\\
      \frac\gamma{x\sigma} & \mbox{for~$x > \frac{\gamma}{\sigma}$}
    \end{cases}
\end{align*}
Applying Lemma~\ref{lem:combLemma} while observing
\begin{itemize}
\item $g' (t) = 0$ for $t \in \left[0,\frac{\gamma}{\sigma}\right]$, and
\item $f (x/t) = 1$ for $t \geq x\sigma/\alpha $,
\end{itemize}
we obtain 
\begin{align*}
  \prob{}{A B \geq x}
    &\leq \int_{\gamma/\sigma}^{x\sigma/\alpha}
      \frac{\alpha + \beta \sqrt{\ln (x\sigma/t)}}{x\sigma/t}
      \frac{\gamma}{t^2\sigma}\diff{t}
      + \int_{x\sigma/\alpha}^\infty \frac{\gamma}{\sigma t^2}\diff{t} \\
    &= \int_{\gamma/\sigma}^{x\sigma/\alpha}
      \frac{\alpha + \beta \sqrt{\ln (x\sigma/t)}}{x\sigma^2}
      \frac{\gamma}{t}\diff{t}
      + \frac{\alpha\gamma}{x\sigma^2}\\
  \intertext{(substituting $s = \sqrt{\ln (x\sigma/t)}, t = x\sigma e^{-s^2}$, 
 which is defined as $x \geq \gamma /\sigma^{2}$, )}
    &= \int_{\sqrt{\ln (x\sigma^2\!/\gamma)}}^{\sqrt{\ln\alpha}}
      \frac{\alpha + \beta s}{x\sigma^2}\frac{\gamma}{x\sigma e^{-s^2}}
      x\sigma (-2s e^{-s^2})\diff{s}
      + \frac{\alpha\gamma}{x\sigma^2}\\
    &= \frac{\gamma}{x\sigma^2}
      \int_{\sqrt{\ln\alpha}}^{\sqrt{\ln (x\sigma^2\!/\gamma)}}
        2s(\alpha + \beta s)\diff{s}
      + \frac{\alpha\gamma}{x\sigma^2}\\
    &= \frac{\alpha\gamma}{x\sigma^2}
      \left(1 + \ln\left(\frac{x\sigma^2}{\alpha\gamma}\right)
        + \frac{2\beta}{3\alpha}
	  \left(\ln^{3/2}\left(\frac{x\sigma^2}{\gamma}\right)
	  - \ln^{3/2} \alpha\right) \right)\\
    &\leq  \frac{\alpha\gamma}{x\sigma^2}
      \left(1 
        + \left(\frac{2\beta}{3\alpha} + 1 \right)
	  ln^{3/2}\left(\frac{x\sigma^2}{\gamma} \right)
          \right),
\end{align*}
as $\alpha \geq 1$.
\end{proof}

\begin{lemma}[linear-bounded expectation]\label{lem:linear-expect}
Let $A$, $B$ and $C$ be positive random variables such that
\[
\prob{}{A \geq x }\leq \frac{\alpha}{x},
\]
for some $\alpha >0$, and
\[
  \forall A, \ \prob{}{B \geq x | A} \leq \prob{}{C \geq x}.
\]
Then,
\[
  \prob{}{AB \geq x}\leq \frac{\alpha }{x}\expec{}{C}.
\]
\end{lemma}
\begin{proof}
Let $g (x)$ be the distribution function of $C$.
By Lemma~\ref{lem:combLemma}, we have
\begin{align*}
 \prob{}{AB \ge x} 
& \le \int_0^\infty
      \left(\frac{\alpha t}{x}\right) (- (1-g)'(t))\, \diff{t}\\
& =
  \frac{\alpha }{x}
   \int_0^\infty t(g'(t))\, \diff{t}\\
& =
  \frac{\alpha }{x}
   \expec{}{C}.
\end{align*}


\end{proof}

\begin{corollary}[linear-chi]\label{cor:chiComb}
Let~$A$ a be positive random variable such that
 
\[
 \prob{}{A \ge x} \le \frac\alpha x. 
\]
  for some~$\alpha > 0$.
Let
  $\bb $ be a $d$-dimensional Gaussian random vector
  (possibly depending upon $A$)
  of variance at most $\sigma^{2}$ centered at a vector
  of norm at most $t$, and let
  $B = \pnorm{2}{\bb}$.
Then,
  \[ \prob{}{AB \ge x} 
\leq 
 \frac{\alpha \sqrt{\sigma^{2} d + t^{2}}}{x}
\]
\end{corollary}
\begin{proof}
As $\expec{}{B} \leq \sqrt{\expec{}{B^{2}}}$, 
  and 
  it is known~\cite[p. 277]{StatsEncyc6} that the expected
    value of $B^{2}$---the non-central $\chi^{2}$-distribution
   with non-centrality parameter $\pnorm{2}{\orig{\bb}}^{2}$---is
    $\sigma^{2}d + \pnorm{2}{\orig{\bb}}^{2}$,
the corollary follows from  Lemma~\ref{lem:linear-expect}.
\end{proof}

\begin{lemma}[Linear to log]\label{lem:linear-log}
Let $A$ be a a positive random variable. 
If there exists an
  $A_{0}\geq 1$ and an $\alpha \geq 1$ such that for all $x\geq A_{0}$,
\[
\prob{A}{A\geq x} \leq \frac{\alpha}{x}.
\]
Then, 
\[
\expec{A}{\max(0,\ln A)} \leq \ln \max(A_{0},\alpha) + 1.
\]
\end{lemma}
\begin{proof}
\begin{eqnarray*}
\expec{A}{\max(0,\ln A)} & = &  \int_{x=0}^{\infty}\prob{A}{\max(0,\ln
A)\geq x} \diff{x}  \\
& \leq  &
\int_{x=0}^{\ln \max(A_{0},\alpha)}1 \diff{x} + 
\int_{x=\ln \max(A_{0},\alpha)}^{\infty}\prob{A}{\ln A\geq x} dx \\
& \leq  & \int_{x=0}^{\ln \max(A_{0},\alpha) } \diff{x} + 
\int_{x=\ln \max(A_{0},\alpha) }^{\infty} \alpha e^{-x} dx \\
& \leq  &
  \ln \max(A_0,\alpha) + 1.
\end{eqnarray*}
\end{proof}

\end{document}